\documentclass[journal,10pt]{IEEEtran}
\usepackage[named]{algo}
\usepackage{graphicx}
\usepackage{amsmath}
\usepackage{amsfonts}
\usepackage{amssymb}
\usepackage{epsfig}
\usepackage{url}

\newcommand{\bo}[1]{{\bf #1}}

\newcommand{\bPsi}{\mbox{\boldmath $\Psi$}}

\newcommand{\bzeta}{\mbox{\boldmath $\zeta$}}

\newcommand{\bm}[1]{\mbox{\boldmath $#1$}}

\newcommand{\bgamma}{\mbox{\boldmath $\gamma$}}

\newcommand{\bnu}{\mbox{\boldmath  $\nu$}}

\newtheorem{theorem}{Theorem}
\newtheorem{proposition}{Proposition}

\newtheorem{corollary}{Corollary}

\newtheorem{proof}{Proof}

\title{Restoration of Poissonian Images Using \\ Alternating
Direction Optimization}

\author{M\'ario A. T. Figueiredo,~\IEEEmembership{Fellow,~IEEE,} \hspace{0.5cm}
Jos\'{e} M. Bioucas-Dias,~\IEEEmembership{Member,~IEEE},
\thanks{M. Figueiredo  and J. Bioucas-Dias are with the {\it Instituto de Telecomunica\c{c}\~{o}es},
   {\it Instituto Superior T\'ecnico,} 1049-001 Lisboa, Portugal. Email:
   {\tt mario.figueiredo@lx.it.pt} and {\tt bioucas@lx.it.pt}.}
\thanks{This work was partially supported by {\it Funda\c{c}\~{a}o
para a Ci\^{e}ncia e Tecnologia} (FCT), Portuguese Ministry of Science
and Higher Education, under projects POSC/EEA-CPS/61271/2004
and PTDC/EEA-TEL/104515/2008.}
   }

\begin{document}

\maketitle

\begin{abstract}
Much research has been devoted to the problem of restoring
Poissonian images, namely for medical and astronomical applications.
However, the restoration of these images using state-of-the-art
regularizers (such as those based on multiscale representations or total
variation) is still an active research area, since the associated optimization
problems are quite challenging. In this paper, we propose an
approach to deconvolving Poissonian images, which is based on an alternating
direction optimization method. The standard regularization (or maximum a posteriori)
restoration criterion, which combines the Poisson log-likelihood with a
(non-smooth) convex regularizer (log-prior), leads to hard
optimization problems: the log-likelihood is non-quadratic and
non-separable, the regularizer is non-smooth, and there is a
non-negativity constraint. Using standard convex analysis tools,
we present sufficient conditions for existence and uniqueness of
solutions of these optimization problems, for several types of
regularizers: total-variation, frame-based analysis, and frame-
based synthesis. We attack these problems with an instance
of the alternating direction method of multipliers (ADMM), which
belongs to the family of augmented Lagrangian algorithms. We study
sufficient conditions for convergence and show that these are
satisfied, either under total-variation or frame-based (analysis and
synthesis) regularization.
The resulting algorithms are shown to outperform alternative
state-of-the-art methods, both in terms of speed and restoration
accuracy.
\end{abstract}


\section{Introduction}
A large fraction of (if not all) the work on image denoising, restoration, and
reconstruction has been devoted to developing regularizers (priors,
from a Bayesian point of view) to deal with the presence of noise
and/or the ill-conditioned or ill-posed nature of the observation
operator, and to devising
efficient algorithms to solve the resulting optimization problems.
Much of that work assumes linearity of the observation
operator ({\it e.g.}, the convolution with some point spread function,
the acquisition of tomographic projections, or simply an identity in
the case of denoising) and the presence of
additive Gaussian noise. For this classical scenario,
recent state-of-the-art methods adopt non-smooth convex
regularizers, such as total-variation or the $\ell_1$ norm of frame
coefficients; the resulting optimization problems are convex but
non-smooth, of very high dimensionality, and have stimulated a considerable
amount of research on special purpose algorithms (see
\cite{FISTA,TwIST,CombettesSIAM,FigueiredoNowak2003,Starck2003b,SpaRSA_SP}
and the many references therein).

The algorithms developed for linear operators and Gaussian noise
cannot be directly applied to other observation models,
such as the Poisson case considered in this paper. Poissonian image
models are well studied and highly relevant in fields such as astronomical
\cite{StarckBook}, biomedical \cite{Dey,FesslerHero,Nowak2000,Sarder,Timmermann,Willett},  and photographic
imaging \cite{Foi_ICIP2005}. A very recent and comprehensive overview of
deconvolution methods for Poissonian images can be found in \cite{DupeFadiliStarck}
(where a state-of-the-art algorithm is also introduced); we refer the
reader to that publication for more references on this topic.

The standard criterion for deconvolution of Poissonian images
consist of a convex constrained optimization problem: the
objective function includes the so-called data term, which is
convex and smooth, but not quadratic, plus a convex non-smooth regularizer
(the log-likelihood and log-prior, from a Bayesian inference perspective),
and a constraint forcing the solution to be non-negative.
Although the problem  is convex, its very high dimensionality (when
dealing with images) usually rules out the direct application of
off-the-shelf optimization algorithms.

Furthermore, the Poisson log-likelihood, which is non-quadratic and
non-separable (except in the pure denoising case) raises several
difficulties to the current state-of-the-art image deconvolution
algorithms. More specifically, the Poisson log-likelihood does
not have a Lipschitz-continuous
gradient, a sufficient condition for the applicability
(with guaranteed convergence) of algorithms of the forward-backward splitting (FBS)
family \cite{CombettesSIAM,DupeFadiliStarck,SpaRSA_SP}. If, nevertheless,
an FBS method is applied, it is well known to be slow, specially
if the observation operator is  severely ill-conditioned, a fact which has stimulated recent
research on faster methods \cite{FISTA,TwIST,SpaRSA_SP}; these
faster algorithms also require the log-likelihood to have a
Lipschitz-continuous gradient, which is not the case with Poissonian
observations.

In this paper, we propose a new approach to tackle the
optimization problem referred to in the previous paragraph.
The proposed algorithm is based on an instance of the
{\it alternating direction method of multipliers}
(ADMM)  \cite{EcksteinBertsekas}, \cite{Gabay}, \cite{Glowinski},
which belongs to the family of augmented Lagrangian methods \cite{Nocedal}.
For this reason, we call it PIDAL (Poisson image deconvolution
by augmented Lagrangian). Although the proposed approach is related
to the recent split-Bregman (SB) technique \cite{GoldsteinOsher},
our splitting strategy and resulting algorithm are quite different
from the one in \cite{GoldsteinOsher} (which, moreover, is not adequate for
Poissonian image models). Finally, we mention that this paper
is an extension of our much shorter and preliminary work \cite{FigueiredoDiasSSP09}.

In recent work, Douglas-Rachford splitting (DRS) methods were proposed to
attack problems in which log-likelihood the does not have a Lipschitz-continuous
gradient  \cite{CombettesPesquet}. In fact, the ADMM is closely related to
DRS methods \cite{EcksteinBertsekas}, \cite{Esser}, so the method proposed
in this paper can also be interpreted from a DRS viewpoint.

In this paper, we will consider three types of regularization:
total variation \cite{Chambolle04,art:Rdin:O:F:Physica:92} and
both frame-based  analysis and frame-based synthesis formulations
 \cite{EladMilanfarRubinstein}.
In Section \ref{sec:problem_formulation}, after
presenting these three formulations, we
derive sufficient conditions for existence and uniqueness of
solutions of the corresponding optimization problems.
The ADMM framework is reviewed in Section~\ref{sec:ADMM},
where we also introduce the particular variant that is
suitable for a linear combination of several convex functions,
which is the form of the objective
function in hand. In Sections~\ref{sec:PIDAL-TV} and \ref{sec:PIDAL-F},
we instantiate the proposed variant of ADMM to the three types
of regularizers considered and provide sufficient conditions
for convergence. Finally, the effectiveness of the resulting
algorithm is illustrated in comparison with current state-of-the-art
alternatives \cite{Dey}, \cite{DupeFadiliStarck}, \cite{Foi_ICIP2005},
\cite{SetzerSteidlTeuber}, \cite{Starck1995}, via a set of
experiments reported in Section~\ref{sec:experiments}.

\section{Problem Formulation}
\label{sec:problem_formulation}
In this section, we begin by reviewing the derivation of the
standard log-likelihood resulting from assuming that the observations
are Poisson distributed with a mean intensity linearly
related with the underlying image to be estimated.
Then, we present three different regularization/Bayesian criteria, using synthesis
and analysis formulations \cite{EladMilanfarRubinstein}, and
study existence/uniqueness of the corresponding solutions.

\subsection{The Linear/Poisson Observation Model}
Let $\bo{y} = (y_1,...,y_m) \in \mathbb{N}_{0}^{m}$ denote an $m$-vector of
observed counts ($\mathbb{N}_0 = \mathbb{N} \cup \{0\}$), assumed to be a sample
of a random vector ${\bf Y} = (Y_1,...,Y_m)$ of $m$ independent Poisson variables,
with probability distribution
\begin{equation}
P[\bo{Y} = \bo{y} | \bm{\lambda} ] = \prod_{i=1}^m \frac{\lambda_i^{y_i} \, e^{-\lambda_i}}{y_i !},
\label{Poisson}
\end{equation}
where $\bm{\lambda} = (\lambda_1,...,\lambda_m) \in \mathbb{R}_+^m$ ($\mathbb{R}_+$
denotes the non-negative reals) is the underlying mean
(intensity) vector,
assumed to be a linear observation of an unknown image ${\bf x}$, {\it i.e.},
\begin{equation}
\bm{\lambda} = {\bf K\, x},\label{blur}
\end{equation}
where ${\bf K}$ the observation operator, which in our finite dimensional
setting is simply a matrix ${\bf K}\in \mathbb{R}^{m\times n}$.
This matrix may model a convolution or some other linear observation
mechanism, such as emission tomography. So that the underlying
unknown ${\bf x}$ can also have the meaning of intensity, it is
commonly assumed that ${\bf x} \in \mathbb{R}_+^n$.
It is usually further assumed that all the elements of
${\bf K}$ are non-negative \cite{Dey,FesslerHero,StarckBook}.
When dealing with images, we adopt the usual vector notation obtained
by stacking the pixels into a vector, in lexicographic order.

Combining (\ref{Poisson}) and (\ref{blur}) and taking logarithms leads to
the negative log-likelihood function \cite{Dey,StarckBook},
\begin{eqnarray}
- \log P[\bo{Y} = \bo{y} | {\bf x} ] &\!\!\! = &\!\!\! \sum_{i=1}^m  ({\bf K\, x})_i - y_i
\log\left( ({\bf K\, x})_i\right) + \log(y_i!)\nonumber \\
&\!\!\! = &\!\!\! {\cal L}({\bf K\, x}) \label{eq:log_likelihood}\\
&\!\!\! = &\!\!\! ( {\cal L}\circ{\bf K}) \, ({\bf x}),
\end{eqnarray}
where $({\bf v})_i$ (or $v_i$) denotes the $i$-th component of some vector ${\bf v}$
and ${\cal L}:\mathbb{R}^m \rightarrow \bar{\mathbb{R}} = \mathbb{R}\cup
\{-\infty,+\infty\}$ is the negative log-likelihood function for the
case ${\bf K=I}$, that is
\begin{equation}
{\cal L}({\bf z}) = \sum_{i=1}^m  -\, y_i\, \log\left( z_i\right) + z_i + \log(y_i!).
\end{equation}

Dealing with the particular case $z_i = 0$ requires some care, because
of the presence of the logarithm. Seen as  function of ${\bf z}$
to be used in a minimization problem, it is convenient to write the
negative log-likelihood function as
\begin{equation}
{\cal L}({\bf z}) = C + \sum_{i=1}^m \xi (z_i , y_i ),
\end{equation}
where $C = \sum_i \log(y_i!)$ is a finite (recall that $0! = 1$)
irrelevant (independent of ${\bf z}$) constant
and $\xi:\mathbb{R} \times \mathbb{N}_0 \rightarrow
\bar{\mathbb{R}}$ is defined as
\begin{equation}
\xi(z,y) = z + \iota_{\,\mathbb{R}_+}(z) - y\, \log(z_+),
\end{equation}
where $\iota_S$ is the indicator function of set $S$,
\[
\iota_S(z) = \left\{\begin{array}{lll} 0 & \Leftarrow & z \in S\\
+\infty & \Leftarrow & z \not\in S, \end{array} \right.
\]
$z_+ = \max\{0,z\}$,  $\log(0)=-\infty$, and $0\log(0) = 0$.

The following two propositions characterize
$\xi$ as a function of its first argument, as well as ${\cal L}$
and ${\cal L}\circ {\bf K}$, in terms of the key
concepts of convex analysis (see Appendix A).

\vspace{0.2cm}
\begin{proposition}
{\sl For any $y\in \mathbb{N}_0$ the function
$\xi(\cdot,y):\mathbb{R}\rightarrow \bar{\mathbb{R}}$ is proper,
lower semi-continuous (lsc), coercive, and convex. If $y>0$, then
 $\xi(\cdot,y)$ is also strictly convex.}
\end{proposition}

\vspace{0.2cm}
\begin{proof}
For $y=0$,  $\xi(z,0) = z + \iota_{\,\mathbb{R}_+}(z)$,
thus $\xi(\cdot,0)$ is the sum of the identity function with
$\iota_{\,\mathbb{R}_+}$, which are both proper,
lsc, coercive, and convex.
For any $y>0$, $\xi(z,y) = z + \iota_{\,\mathbb{R}_+}(z) - y\log (z_+)$; since
$y\log ((\cdot)_+)$ is also proper, lsc, coercive, and convex,
so is $\xi(\cdot,y)$. Finally, if $y>0$, $y\log ((\cdot)_+)$ is
strictly convex (see the definition in Appendix A), thus
$\xi(\cdot,y)$ is also strictly convex.
\end{proof}

\vspace{0.2cm}
\begin{proposition}
\label{prop:about_L}
{\sl Function ${\cal L}$ is proper, lsc, coercive, and convex. If  $y_i\neq 0$, for $i=1,...,m$,
${\cal L}$ is also strictly convex. Function ${\cal L}\circ{\bf K}$
is proper, lsc, and convex. Function ${\cal L}\circ{\bf K}$ is
coercive if ${\bf K}$ is injective. Function ${\cal L}\circ{\bf K}$ is
strictly convex if ${\bf K}$ is injective and $y_i\neq 0$, for $i=1,...,m$.}
\end{proposition}

\vspace{0.2cm}
\begin{proof}
Function ${\cal L}$ is the sum of proper, lsc, coercive, convex functions.
If  $y_i\neq 0$, for $i=1,...,m$,  the functions in the sum are
also strictly convex, thus ${\cal L}$ is also strictly convex.
Function ${\cal L}\circ{\bf K}$ is the composition of a proper, lsc, convex functions
with a linear function, thus it is proper, lsc, and convex. If ${\bf K}$ is
injective, its null set is the zero vector, thus
$\lim_{\|{\bf x}\|\rightarrow +\infty} \|{\bf K\, x}\| = +\infty$, thus
${\cal L}\circ{\bf K}$ is coercive. Finally, if ${\bf K}$ is
injective and $y_i\neq 0$, for $i=1,...,m$, ${\cal L}$ is strictly convex,
thus so is ${\cal L}\circ{\bf K}$.
\end{proof}

\subsection{Regularization Criteria: Analysis and Synthesis Formulations}
Under a {\it maximum a posteriori} (MAP) or regularization criterion, the
image estimate is obtained by solving a variational problem: minimizing an
objective function, which includes the log-likelihood term plus
a regularizer \cite{FigueiredoDiasSSP09}, \cite{LeChartrand}, \cite{StarckBook},
under a positivity constraint. We will
now describe three possible ways of building such an
objective function.

\vspace{0.2cm}
\subsubsection{Total Variation Regularization}
A standard choice for regularization of digital image restoration/reconstruction
problems is the isotropic discrete TV regularizer \cite{Chambolle04},
\begin{equation}
\mbox{TV}(\bo{x})  =  \sum_{s=1}^n \sqrt{(\Delta^h_s\bo{x})^2+(\Delta^v_s\bo{x})^2},\label{eq:theTV}
\end{equation}
where  $(\Delta^h_s\bo{x}$ and $\Delta^v_s\bo{x})$ denote the horizontal and vertical first
order differences at pixel $s$, respectively. This regularizer is a discrete
version of the TV regularizer proposed in \cite{art:Rdin:O:F:Physica:92}.
The  resulting optimization problem is
\begin{equation}
 \min_{{\bf x}} \, L^{\mbox{\scriptsize TV}}({\bf x}) \label{eq:L_uncons_TV}
\end{equation}
with
\begin{equation}
L^{\mbox{\scriptsize TV}}({\bf x}) = {\cal L}(\bo{K\, x}) + \tau \mbox{ TV}({\bf x}) + \iota_{\mathbb{R}_+^n}({\bf x}),
\label{eq:objective_TV}
\end{equation}
where  $\tau \in \mathbb{R}_+$ is the regularization parameter and the role
of $\iota_{\mathbb{R}_+^n}$, the indicator of the first orthant, is to
impose the non-negativity constraint on the estimate. The next proposition
concerns the existence and uniqueness of minimizers of $L^{\mbox{\scriptsize TV}}$.

\vspace{0.2cm}
\begin{proposition}
Consider the function $L^{\mbox{\scriptsize TV}}$ defined in (\ref{eq:objective_TV}).
\begin{description}
\item[a)] $L^{\mbox{\scriptsize TV}}$ is proper, lsc, and convex.
\item[b)] If the intersection of the space of constant images
$\{ {\bf x} = \alpha(1,1,...,1),\;  \alpha\in\mathbb{R}\}$ (which is
the null space of $\mbox{TV}$) with the null space
of ${\bf K}$ is just the zero vector, then $L^{\mbox{\scriptsize TV}}$ is coercive,
and (\ref{eq:L_uncons_TV}) has a solution.
\item[c)] If ${\bf K}$ is injective, then (\ref{eq:L_uncons_TV}) has a solution.
\item[d)] If ${\bf K}\in\mathbb{R}_+^{m\times n}$, and at least one element of
${\bf K}$ is strictly positive, then (\ref{eq:L_uncons_TV}) has a solution.
\item[e)] If ${\bf K}$ is injective and $y_i\neq 0$, for $i=1,...,m$, then
$L^{\mbox{\scriptsize TV}}$ is coercive and strictly convex thus there is a unique solution.
\end{description}
\end{proposition}

\begin{proof}
\begin{description}
\item[a)] The  functions $\iota_{\mathbb{R}_+^n}$,  $\mbox{TV}$,
and ${\cal L}\circ{\bf K}$ (Proposition \ref{prop:about_L}) are proper, lsc, and convex,
thus so is their sum.
\item[b)] Similar to  \cite{ChambolleLions}.
\item[c)] If ${\bf K}$ is injective, its null space is just the zero vector, thus
${\cal L}\circ{\bf K}$ and $L^{\mbox{\scriptsize TV}}$ are coercive.
\item[d)] If all the elements of ${\bf K}$ are non-negative and at least one is
positive, then the constant vector $(1,1,...,1)$ doesn't belong to the null
space of ${\bf K}$ and the result follows from (b).
\item[e)] If ${\bf K}$ is injective and  $y_i\neq 0$, for $i=1,...,m$,
${\cal L}\circ{\bf K}$ is strictly convex (Proposition~\ref{prop:about_L}), thus
so is $L^{\mbox{\scriptsize TV}}$ and its minimizer is unique.
\end{description}
\end{proof}

\subsubsection{Frame Analysis Regularization}
The use of a regularizer which is a direct function of the unknown image
(as in  (\ref{eq:L_uncons_TV})--(\ref{eq:objective_TV})) corresponds to a so-called analysis-based
prior/regularizer \cite{EladMilanfarRubinstein}.  Another well-known type of analysis-based
regularization penalizes the norm (typically $\ell_1$)
of the representation coefficients of ${\bf x}$ on some wavelet basis or
tight frame, given by ${\bf P x}$, where ${\bf P}$ is the analysis operator
associated with the frame \cite{Mallat}.
This approach leads to the following optimization problem:
\begin{equation}
       \min_{\bo{x}} \, L^{\mbox{\scriptsize FA}}({\bf x}), \label{eq:L_uncons_analysis}
\end{equation}
where FA stands for {\it frame analysis} and
\begin{equation}
L^{\mbox{\scriptsize FA}}({\bf x}) = {\cal L}(\bo{K\, x}) + \tau \|{\bf P} {\bf x}\|_1 + \iota_{\mathbb{R}_+^n}({\bf x});\label{eq:frame_analysis}
\end{equation}
as above, $\tau$ is the regularization parameter and
$\iota_{\mathbb{R}_+^n}$ imposes the non-negativity constraint on the estimate.
The next proposition addresses the existence and uniqueness of minimizers of $L^{\mbox{\scriptsize FA}}$.

\vspace{0.2cm}
\begin{proposition}\label{prop:existenceFA}
Consider the function $L^{\mbox{\scriptsize FA}}$ defined in (\ref{eq:frame_analysis}).
\begin{description}
\item[a)] $L^{\mbox{\scriptsize FA}}$ is proper, lsc, convex, and coercive, thus
has a minimizer.
\item[b)] If ${\bf K}$ is injective and $y_i\neq 0$, for $i=1,...,m$, then
$L^{\mbox{\scriptsize FA}}$ is strictly convex with a unique minimizer.
\end{description}
\end{proposition}

\begin{proof}
\begin{description}
\item[a)] The  functions $\iota_{\mathbb{R}_+^n}$,  $\|\cdot \|_1\circ {\bf P}$,
and ${\cal L}\circ{\bf K}$ (Proposition \ref{prop:about_L}) are proper, lsc, and convex,
thus so is their sum. Furthermore, since ${\bf P}$ is the analysis
operator of a tight frame, its null space is simply the zero vector, thus
$\|\cdot \|_1\circ {\bf P}$ is coercive.
\item[b)] If ${\bf K}$ is injective and  $y_i\neq 0$, for $i=1,...,m$,
${\cal L}\circ{\bf K}$ is strictly convex (Proposition~\ref{prop:about_L}), thus
so is $L^{\mbox{\scriptsize FA}}$ and its minimizer is unique.
\end{description}
\end{proof}

\vspace{0.2cm}
\subsubsection{Frame Synthesis Regularization}
Finally, another well-known class of approaches is known as synthesis-based
\cite{EladMilanfarRubinstein}. Here, the unknown image is represented
on a frame ({\it e.g.}, of wavelets, curvelets, or other multi-scale
system) and then the coefficients of this representation are estimated
from the observed data, under some regularizer. With ${\bf W}\in \mathbb{R}^{n\times d}$
denoting the synthesis matrix of the frame, the image is written
as ${\bf x = W\, s}$, where ${\bf s}$ is the vector of representation
coefficients, and the resulting optimization problem is
\begin{equation}
       \min_{\bo{s}} \,  L^{\mbox{\scriptsize FS}}({\bf s}) \label{eq:L_uncons_synthesis}
\end{equation}
where FS stands for {\it frame synthesis} and
\begin{equation}
\label{eq:FS_objective}
L^{\mbox{\scriptsize FS}}({\bf s}) = {\cal L}(\bo{K\, W\, s}) + \tau \|{\bf s}\|_1 + \iota_{\mathbb{R}_+^n}({\bf W\, s}).
\end{equation}
Naturally, the indicator function $\iota_{\mathbb{R}_+^n}$ forcing the image estimate
to be non-negative is applied to the image ${\bf W\, s}$ and not its coefficients.
The next proposition addresses the existence and uniqueness of minimizers of $L^{\mbox{\scriptsize FS}}$.

\vspace{0.2cm}
\begin{proposition}
Consider the function $L^{\mbox{\scriptsize FS}}$ defined in (\ref{eq:FS_objective}).
\begin{description}
\item[a)] $L^{\mbox{\scriptsize FS}}$ is proper, lsc, convex, and coercive, thus
has a minimizer.
\item[b)] If ${\bf K\, W}$ is injective and $y_i\neq 0$, for $i=1,...,m$, then
$L^{\mbox{\scriptsize FS}}$ is strictly convex with a unique minimizer.
\end{description}
\end{proposition}

\begin{proof}
\begin{description}
\item[a)] The  functions $\iota_{\mathbb{R}_+^n}\circ {\bf W}$,  $\|\cdot \|_1$,
and ${\cal L}\circ{\bf K\,W}$ (Proposition \ref{prop:about_L}) are proper, lsc, and convex,
thus so is their sum. Furthermore, since $\|\cdot \|_1$ is coercive,
$L^{\mbox{\scriptsize FS}}$ is coercive.
\item[b)] Same proof as that of Proposition~\ref{prop:existenceFA} (b).
\end{description}
\end{proof}

\section{The Alternating Direction Method of Multipliers (ADMM)}
\label{sec:ADMM}
\subsection{The Standard ADMM}
The key tool in this paper is the {\it alternating direction method of multipliers}
(ADMM)  \cite{EcksteinBertsekas}, \cite{Gabay}, \cite{Glowinski}.
Consider an unconstrained problem of the form
\begin{equation}
\min_{{\bf z}\in \mathbb{R}^d} \; f_1({\bf z}) + f_2 ({\bf G\, z}),\label{unconstrained_basic}
\end{equation}
where $f_1:\mathbb{R}^d \rightarrow \bar{\mathbb{R}}$, $f_2:\mathbb{R}^p
\rightarrow \bar{\mathbb{R}}$, and ${\bf G}\in\mathbb{R}^{p\times d}$.
The ADMM for this problem is defined in Fig.~\ref{fig:ADMM}.

\begin{figure}[hbt]
\begin{center}
\fbox{\parbox{0.9\columnwidth}{\small
\begin{algorithm}{ADMM}{
\label{alg:salsa2}}
Set $k=0$, choose $\mu > 0$, ${\bf u}_0$,  and  ${\bf d}_0$.\\
\qrepeat\\
   $  {\bf z}_{k+1}  \in  \arg\min_{{\bf z}} f_{1}({\bf z})
 + \frac{\mu}{2} \|{\bf G\, z} - {\bf u}_k - {\bf d}_k\|_2^2$\rule[-0.25cm]{0cm}{0.5cm}\\
  $  {\bf u}_{k+1}  \in  \arg\min_{{\bf u}} f_{2}({\bf u})
 + \frac{\mu}{2} \|{\bf G\, z}_{k+1} - {\bf u} - {\bf d}_k\|_2^2$\rule[-0.25cm]{0cm}{0.5cm}\\
     ${\bf d}_{k+1} \leftarrow {\bf d}_{k} - ({\bf G \, z}_{k+1}  - {\bf u}_{k+1})$\rule[-0.25cm]{0cm}{0.5cm}\\
     $k \leftarrow k+1$
\quntil stopping criterion is satisfied.
\end{algorithm} }}
 \caption{The ADMM algorithm.}\label{fig:ADMM}
\end{center}
\end{figure}

For later reference, we now recall a theorem by Eckstein and Bertsekas
in which convergence of (a generalized version of) ADMM is shown.
\vspace{0.3cm}
\begin{theorem}[Eckstein-Bertsekas, \cite{EcksteinBertsekas}]
\label{th:Eckstein}{\sl Consider problem  (\ref{unconstrained_basic}), where
 ${\bf G}\in\mathbb{R}^{p\times d}$ has full column rank and
 $f_1\,:\,\mathbb{R}^d\rightarrow \bar{\mathbb{R}}$ and $f_2\,:\,\mathbb{R}^p\rightarrow
\bar{\mathbb{R}}$ are closed, proper, convex functions. Consider arbitrary
 $\mu>0$ and ${\bf u}_0, {\bf d}_0\in \mathbb{R}^p$.
 Let $\{\eta_k \geq 0, \; k=0,1,...\}$ and $\{\rho_k \geq 0,
\; k=0,1,...\}$ be two sequences such that
\[
\sum_{k=0}^\infty \eta_k < \infty \;\;\;\mbox{and} \;\;\; \sum_{k=0}^\infty \rho_k <
\infty.
\]
Consider three sequences $\{{\bf z}_k \in \mathbb{R}^{d}, \; k=0,1,...\}$, $\{{\bf u}_k
\in \mathbb{R}^{p}, \; k=0,1,...\}$, and $\{{\bf d}_k \in \mathbb{R}^{p}, \; k=0,1,...\}$
that satisfy
\begin{eqnarray}
  \left\| {\bf z}_{k+1} - \arg\min_{{\bf z}} f_{1}({\bf z})
 + \frac{\mu}{2} \|{\bf G}{\bf z} \! - \!{\bf u}_k \! -\! {\bf d}_k\|_2^2 \right\| & \leq & \eta_k \nonumber\\
  \left\| {\bf u}_{k+1}  - \arg\min_{{\bf u}} f_{2}({\bf u})
 + \frac{\mu}{2} \|{\bf G}{\bf z}_{k+1} \! - \! {\bf u} \! - \! {\bf d}_k\|_2^2 \right\| & \leq &  \rho_k \nonumber
\end{eqnarray}
and
\begin{equation}
 {\bf d}_{k+1}  =  {\bf d}_{k} - ({\bf G\, z}_{k+1} - {\bf u}_{k+1}).
 \end{equation}
Then, if (\ref{unconstrained_basic}) has a solution, say ${\bf z}^*$,
the sequence $\{{\bf z}_k\}$  converges to ${\bf z}^*$.
If (\ref{unconstrained_basic}) does not have a solution, then at least
one of the sequences $\{{\bf u}_k \}$ or $\{{\bf d}_k\}$ diverges.}
\end{theorem}
\vspace{0.3cm}

According to Theorem~\ref{th:Eckstein}, it is not necessary to
exactly solve the minimizations in lines 3 and 4 of ADMM: as long as
the sequences of errors are absolutely summable, convergence is not
compromised. As shown in Section~\ref{sec:PIDAL-TV-conv-approx},
this fact is quite relevant in designing instances of
ADMM, when these minimizations lack closed form solutions.

The proof of Theorem \ref{th:Eckstein} is based on the equivalence between ADMM and the
DRS method applied to the dual of  problem
(\ref{unconstrained_basic}).  For recent and comprehensive reviews of ADMM, DRS,
and their relationship with Bregman and split-Bregman methods, see \cite{Esser},
\cite{Setzer}.

\subsection{A Variant of ADMM}\label{sec:variant}
Notice that the ADMM and the associated convergence theorem presented
in the previous subsection apply to objective functions of the form
(\ref{unconstrained_basic}), {\it i.e.}, which are the sum of two
functions. The fact that our objective functions (\ref{eq:L_uncons_TV}),
(\ref{eq:L_uncons_analysis}), and (\ref{eq:L_uncons_synthesis}) involve
more than two terms raises the following question: what is the best
way of mapping an objective with more than two terms into (\ref{unconstrained_basic})
so that the resulting ADMM is easily applicable and the conditions
of Theorem \ref{th:Eckstein} still hold. In this section,
we give an answer to this question, which will constitute the core of
our approach.

Consider a generalization
of problem (\ref{unconstrained_basic}), where instead of two functions,
we have $J$ functions, that is,
\begin{equation}
\min_{{\bf z}\in \mathbb{R}^d} \; \sum_{j=1}^J g_j ({\bf H}^{(j)}\,{\bf z}),\label{unconstrained_compound}
\end{equation}
where $g_j : \mathbb{R}^{p_j} \rightarrow \bar{\mathbb{R}}$ are closed, proper,
convex functions, and ${\bf H}^{(j)} \in \mathbb{R}^{p_j \times d}$ are arbitrary
matrices. The minimization problem (\ref{unconstrained_compound}) can be written
as (\ref{unconstrained_basic}) using the following correspondences: $f_1 = 0$,
\begin{equation}
{\bf G}  =  \left[ \begin{array}{c} {\bf H}^{(1)} \\ \vdots \\
{\bf H}^{(J)} \end{array}\right] \in \mathbb{R}^{p\times d},\label{eq:stackedG}
\end{equation}
where $p=p_1 + \dots + p_J$, and $f_2 : \mathbb{R}^{p\times d} \rightarrow \bar{\mathbb{R}}$ given by
\begin{equation}
f_2({\bf u}) = \sum_{j=1}^J g_j ({\bf u}^{(j)}),\label{eq:f_2}
\end{equation}
where ${\bf u}^{(j)} \in \mathbb{R}^{p_j}$ and ${\bf u} = [({\bf u}^{(1)})^T,\dots,({\bf u}^{(J)})^T]^T \in \mathbb{R}^{p}$.

We are now in position to apply ADMM. The resulting algorithm has
exactly the same structure as the one in Fig.~\ref{fig:ADMM}
with
 \[
{\bf d}_k = \left[ \begin{array}{c}{\bf d}^{(1)}_k\\ \vdots \\{\bf d}^{(J)}_k\end{array}\right], \;\;\;
{\bf u}_k = \left[ \begin{array}{c}{\bf u}^{(1)}_k\\ \vdots \\{\bf u}^{(J)}_k\end{array}\right].
\]

The fact that $f_1 = 0$ turns Step 3 of the algorithm into a
simple quadratic minimization problem, which has a unique solution if ${\bf G}$
has full column rank:
\begin{eqnarray}
\arg\min_{{\bf z}} \; \bigl\|{\bf G\, z} - \bzeta_k\bigr\|_2^2 & = & \left({\bf G}^T {\bf G}\right)^{-1} {\bf G}^T \bzeta_k,
\label{eq:quadratic_problem}\\
& & \hspace{-2cm} = \; \; \biggl[ \sum_{j=1}^J
({\bf H}^{(j)})^T{\bf H}^{(j)}\biggr]^{-1} \sum_{j=1}^J \bigl({\bf H}^{(j)}\bigr)^T \bzeta^{(j)}_k,\nonumber
\end{eqnarray}
where $\bzeta_k = {\bf u}_k + {\bf d}_k$ (and, naturally, $\bzeta_k^{(j)} = {\bf u}_k^{(j)} + {\bf d}_k^{(j)}$)
and the second equality results from the particular structure of ${\bf G}$ in (\ref{eq:stackedG}).

Furthermore, our particular way of mapping problem (\ref{unconstrained_compound}) into
problem (\ref{unconstrained_basic}) allows decoupling the minimization in Step 4
of ADMM into a set of $J$ independent ones. In fact,
\begin{equation}
 {\bf u}_{k+1}   \leftarrow  \arg\min_{{\bf u}} f_{2}({\bf u})
 + \frac{\mu}{2} \|{\bf G\, z}_{k+1} - {\bf u} - {\bf d}_k\|_2^2 \label{update_u}
\end{equation}
which can be written as
\begin{eqnarray}
\left[ \begin{array}{l} {\bf u}^{(1)}_{k+1}\\ \vdots \\ {\bf u}^{(J)}_{k+1}
  \end{array}\right]
& \!\! \leftarrow  & \!\! \arg\min_{{\bf u}^{(1)}, \dots , {\bf u}^{(J)}} \; g_1 ({\bf u}^{(1)}) + \dots + g_J ({\bf u}^{(J)}) \nonumber \\
  & &\hspace{-1cm} +
\frac{\mu}{2}\left\| \left[ \begin{array}{c} {\bf H}^{(1)} \\ \vdots \\ {\bf H}^{(J)} \end{array}\right]  {\bf z}_{k+1} -
\left[ \begin{array}{c} {\bf u}^{(1)} \\ \vdots \\ {\bf u}^{(J)}\end{array}\right] - \left[\begin{array}{l} {\bf d}^{(1)}_k \\
\vdots \\ {\bf d}^{(J)}_k\end{array}\right]  \right\|_2^2.\nonumber
\end{eqnarray}
Clearly, the minimizations with respect to ${\bf u}^{(1)},\dots,{\bf u}^{(J)}$
are decoupled, thus can be solved separately, leading to
\begin{equation}
 {\bf u}^{(j)}_{k+1}  \leftarrow  \arg\min_{{\bf v}\in \mathbb{R}^{p_j}} \;  g_j ({\bf v}) +
\frac{\mu}{2} \, \bigl\| {\bf v}  - {\bf s}^{(1)}_k\bigr\|_2^2 , \label{eq:Moreau_1}
\end{equation}
for $j=1,...,J$, where
\[
{\bf s}^{(j)}_k =  {\bf H}^{(j)} {\bf z}_{k+1} - {\bf d}^{(j)}_k.
\]
Eq.~(\ref{eq:Moreau_1}) defines the so-called {\it Moreau proximity
operators} \cite{CombettesSIAM} of $g_1,...,g_J$,
applied to ${\bf s}^{(1)}_k,...,{\bf s}^{(J)}_k$, respectively,
denoted as
\begin{equation}
{\bf u}_{k+1}^{(j)} = \bPsi_{g_j/\mu}( {\bf s}^{(j)}_k ).
\end{equation}

Some comments on the algorithm are relevant. Firstly, being exactly
an ADMM, and since all the functions $g_j$, for $j=1,...,J$, are
closed, proper, and convex, convergence is guaranteed if ${\bf G}$
has full column rank. This full column rank condition, which is also
required for the inverse in (\ref{eq:quadratic_problem}) to exist,
will be studied in the next section for each of the specific problems
considered in this paper.

For some functions, this mapping can be computed exactly in
closed form. For example, if $g_j ({\bf x}) = \|{\bf x}\|_1$,
the corresponding proximity operator $\bPsi_{g_j/\mu}$ is simply a soft threshold,
\begin{equation}
\bPsi_{g_j/\mu}({\bf v}) = \mbox{soft}({\bf v},1/\mu) = \mbox{sign}({\bf v}) \odot \max\{ |{\bf v}| - (1/   \mu) ,0\},
\label{eq:soft}
\end{equation}
where $\mbox{sign}(\cdot)$ denotes the component-wise
application of the sign function, $\odot$ denotes the
component-wise product, $|{\bf v}|$ denotes the vector of
absolute values of the elements of ${\bf v}$, and the maximum
is computed in a component-wise fashion. For other functions,
the corresponding Moreau proximity operator does not have a
simple close form solution and needs to be computed numerically.

%

\section{Poissonian Image Reconstruction with TV-Based Regularization}
\label{sec:PIDAL-TV}
\subsection{Applying ADMM}
In this section, we apply the algorithmic framework presented in
Section \ref{sec:variant} to the total-variation-based
criterion (\ref{eq:L_uncons_TV})--(\ref{eq:objective_TV}). The objective function in
(\ref{eq:objective_TV}) has the form (\ref{unconstrained_compound})
with  $J=3$,
\begin{equation}
g_1 \equiv {\cal L}, \;\;\; g_2 \equiv \tau\; \mbox{TV}, \;\;\; g_3 \equiv \iota_{\mathbb{R}_+^n}
\label{eq:g_functions_TV}
\end{equation}
and
\begin{equation}
{\bf H}^{(1)} \equiv {\bf K}, \;\;\; {\bf H}^{(2)} \equiv {\bf I}, \;\;\; {\bf H}^{(3)} \equiv {\bf I}.
\label{eq:G4TV}
\end{equation}
The resulting ADMM algorithm, which we call PIDAL-TV (Poisson image deconvolution
by augmented Lagrangian -- total variation), is shown in Fig.~\ref{fig:PIDAL}.

\begin{figure}
\begin{center}
\fbox{\parbox{0.9\columnwidth}{\footnotesize
\begin{algorithm}{Poisson Image Deconvolution by AL (PIDAL-TV)}{}
Choose $\bo{u}_0^{(1)}$, $\bo{u}_0^{(2)}$, $\bo{u}_0^{(3)}$, $\bo{d}^{(1)}_{0}$, $\bo{d}^{(2)}_{0}$,
$\bo{d}^{(3)}_{0}$, $\mu$, and $\tau$. Set $k \leftarrow 0$.\\
\qrepeat \\
  $\bzeta_k^{(1)} \leftarrow {\bf u}_k^{(1)} + \bo{d}^{(1)}_{k}$\rule[-0.3cm]{0cm}{0.6cm}\\
  $\bzeta_k^{(2)} \leftarrow {\bf u}_k^{(2)} + \bo{d}^{(2)}_{k}$\rule[-0.3cm]{0cm}{0.6cm}\\
  $\bzeta_k^{(3)} \leftarrow {\bf u}_k^{(3)} + \bo{d}^{(3)}_{k}$\rule[-0.3cm]{0cm}{0.6cm}\\
  $\bgamma_k \leftarrow  {\bf K}^T \bzeta_k^{(1)} + \bzeta_k^{(2)} + \bzeta_k^{(3)} $\rule[-0.3cm]{0cm}{0.6cm}\\
  $\bo{z}_{k+1} \leftarrow \left( {\bf K}^T{\bf K} + 2\, {\bf I}\right)^{-1} \bgamma_k
  \rule[-0.3cm]{0cm}{0.6cm}$\\
        $\bnu_k^{(1)} \leftarrow {\bf K \, z}_{k+1} - \bo{d}^{(1)}_{k}$\rule[-0.3cm]{0cm}{0.6cm}\\
${\displaystyle \bo{u}^{(1)}_{k+1} \leftarrow \arg\min_{\bo{v}} \frac{\mu}{2}\| {\bf v} - \bnu_k^{(1)}\|_2^2 + \sum_{i=1}^m \xi(v_i , y_i)  }$\rule[-0.3cm]{0cm}{0.6cm}\\
   $\bnu_k^{(2)} \leftarrow {\bf z}_{k+1}-\bo{d}^{(2)}_{k}$\rule[-0.3cm]{0cm}{0.6cm}\\
${\displaystyle\bo{u}_{k+1}^{(2)} \leftarrow \arg\min_{\bo{v}} \frac{\mu}{2}\|
\bo{v}-\bnu_k^{(2)} \|^2 + \tau \, \mbox{TV}(\bo{v})}$.\rule[-0.3cm]{0cm}{0.6cm}\\
$\bnu_k^{(3)} \leftarrow {\bf z}_{k+1} - \bo{d}^{(3)}_{k}$\rule[-0.3cm]{0cm}{0.6cm}\\
${\displaystyle\bo{u}_{k+1}^{(3)} \leftarrow \arg\min_{\bo{v}} \frac{\mu}{2}\|
\bo{v}-\bnu_k^{(3)} \|^2 + \iota_{\mathbb{R}_+^n}(\bo{v})}$.\rule[-0.3cm]{0cm}{0.6cm}\\
$\bo{d}^{(1)}_{k+1}  \leftarrow  \bo{d}^{(1)}_{k} - (\bo{K\, z}_{k+1} -\bo{u}_{k+1}^{(1)})$\rule[-0.3cm]{0cm}{0.6cm}\\
$\bo{d}^{(2)}_{k+1}   \leftarrow  \bo{d}^{(2)}_{k} - (\bo{z}_{k+1} -\bo{u}_{k+1}^{(2)})$\rule[-0.3cm]{0cm}{0.6cm}\\
$\bo{d}^{(3)}_{k+1}   \leftarrow \bo{d}^{(3)}_{k} - (\bo{z}_{k+1} -\bo{u}_{k+1}^{(3)})$\rule[-0.3cm]{0cm}{0.6cm}\\
$k \leftarrow k + 1$
\quntil some stopping criterion is satisfied.
 \end{algorithm} }}
 \caption{The PIDAL-TV algorithm.}\label{fig:PIDAL}
\end{center}
\end{figure}

\subsection{Implementation Aspects and Computational Cost of PIDAL-TV}
Notice that line 7 of PIDAL-TV corresponds to (\ref{eq:quadratic_problem})
for the particular form of matrix  ${\bf G}$ in this problem:
${\bf G} = [{\bf K}^T \; {\bf I}\;\, {\bf I}\, ]^T$ (see (\ref{eq:g_functions_TV})
and (\ref{eq:G4TV})),
which is of course of full
column rank. Moreover, if $\bf K$ models a periodic convolution, it is a
block circulant matrix and the inversion in line 7  of the algorithm can
be implemented in $O(n\log n)$ operations, via the FFT algorithm.
Although this is a well-known fact, we include the derivation in the
next paragraph, for the sake of completeness.

Assuming that the convolution is periodic, ${\bf K}$ is
block-circulant  with circulant blocks and can be factorized as
\begin{equation}
{\bf K = U}^H {\bf D U},\label{eq:diagonalization}
\end{equation}
where ${\bf U}$ is the matrix that represents the
2D discrete Fourier transform (DFT), ${\bf U}^H = {\bf U}^{-1} $ is its
inverse (${\bf U}$ is unitary, {\it i.e.}, ${\bf UU}^H = {\bf U}^H {\bf U} = {\bf I}$),
and ${\bf D}$ is a diagonal matrix containing the DFT coefficients of
the convolution operator represented by ${\bf K}$.
Thus (with ${\bf K}^T = {\bf K}^H$, since ${\bf K}$ is a real matrix)
\begin{eqnarray}
\left({\bf K}^T{\bf K} + 2 {\bf I}\right)^{-1} & = & \left({\bf U}^H{\bf D}^*{\bf D}{\bf U} + 2 {\bf U}^H{\bf U}\right)^{-1}\\
& = & {\bf U}^H\left(|{\bf D}|^2 + 2\, {\bf I}\right)^{-1}{\bf U},\label{eq:Wiener2}
\end{eqnarray}
where $(\cdot)^*$ denotes complex conjugate and $|{\bf D}|^2$ the
squared absolute values of the entries of  ${\bf D}$.
Since $|{\bf D}|^2 + 2\, {\bf I}$ is diagonal, its inversion has
$O(n)$ cost. Products by ${\bf U}$  and ${\bf U}^H$ can
be carried out with $O(n\log n)$ cost using the FFT algorithm.

The minimization in line 9 of PIDAL-TV is separable. With respect to
each component, it has the form
\begin{equation}
\min_{v} \left\{ v + \iota_{\,\mathbb{R}_+}(v) - y\, \log(v_+) + \frac{\mu}{2}(v - \nu )^2\right\}.
\label{eq:Poisson_Moreau}
\end{equation}
It is simple to show that the solution of (\ref{eq:Poisson_Moreau}) leads to
\begin{equation}
{u}^{(1)}_{i,k+1} = \frac{1}{2} \left( \nu_{i,k}^{(1)} - \frac{1}{\mu} + \sqrt{\left( \nu_{i,k}^{(1)}-\frac{1}{\mu} \right)^2 + 4\, y_i}\; \right), \label{eq:Poisson_Moreau2}
\end{equation}
where $\nu_{i,k}^{(1)}$ denotes the $i$-th component of $\bnu_k^{(1)}$.
Notice that ${u}^{(1)}_{i,k+1}$ is always a non-negative quantity.

The minimization in line 11 of PIDAL-TV is, by definition,
the Moreau proximity operator $\bPsi_{(\tau/\mu)\mbox{\scriptsize TV}}:\mathbb{R}^n\rightarrow \mathbb{R}^n $
\cite{CombettesSIAM}, which corresponds to applying TV-denoising to
$\bnu_k^{(2)}$. Below, we address in detail the issue of how
to implement this operator and its implications to the convergence of PIDAL-TV.
Suffice it to say here that most TV-denoising algorithms have $O(n)$ cost.

The minimization in line 13 of PIDAL-TV corresponds to the projection
of $\bnu_k^{(3)}$ onto the first orthant, thus
\begin{equation}
{\bf u}_{k+1}^{(3)} = \max \left\{ \bnu_{k}^{(3)} , 0 \right\},\label{eq:orthant}
\end{equation}
where the maximum is to be understood in a component-wise sense; this projection has
of course $O(n)$ cost.

From the observations in the previous paragraphs, the computational costs of
the lines of PIDAL-TV are  the following. Lines 3, 4, 5, 9, 10, 11, 12, 13, 15, and 16
have $O(n)$ cost. Lines 6, 7, 8, and 14 have $O(n\log n)$ cost. Thus the
computational cost of PIDAL-TV scales as $O(n\log n)$.

\subsection{Convergence of PIDAL-TV:  Exact TV}
Convergence of PIDAL-TV is addressed by the following corollary
of Theorem \ref{th:Eckstein}, for the (ideal) case where
$\bPsi_{(\tau/\mu)\mbox{\scriptsize TV}}$ (line 11) is computed exactly.
The minimizations in lines 9 and 13 have the closed-form
solutions given in (\ref{eq:Poisson_Moreau2}) and (\ref{eq:orthant}).

\vspace{0.3cm}
\begin{corollary}\label{cor:PIDAL_TV}
{ \sl
If the minimizations in lines 9, 11, and 13 of PIDAL-TV are solved
exactly, then the algorithm converges to a minimizer of  (\ref{eq:objective_TV}),
provided one exists.}
\end{corollary}

\vspace{0.3cm}
\begin{proof}
PIDAL-TV is an instance of ADMM in Fig.~\ref{fig:ADMM}, where $f_1 = 0$
and $f_2$ has the form (\ref{eq:f_2}), with $J=3$ and the $g_i$ functions  given in
(\ref{eq:g_functions_TV}), which are all closed, proper, and convex.
Function $f_2$ is thus also closed, proper, and convex. Matrix
${\bf G} = [{\bf K}^T \; {\bf I}\;\, {\bf I}\, ]^T$ has full column rank.
The minimization in line 4 of ADMM corresponds to lines 9, 11, and 13 of PIDAL-TV;
if these minimizations are solved exactly, then according
to Theorem \ref{th:Eckstein}, convergence to a minimizer of the objective function,
if one exists, is guaranteed.
\end{proof}
\vspace{0.3cm}

\subsection{Convergence of PIDAL-TV:  Approximate TV}
\label{sec:PIDAL-TV-conv-approx}
As is well known, the TV denoising problem has no closed
form solution, and many
iterative algorithms have been proposed to solve it (see \cite{Chambolle04},
\cite{tech:Chan:EPY:UCLA:05}, \cite{Jensen}, \cite{art:Rdin:O:F:Physica:92},
and references therein). Here, we adopt Chambolle's algorithm \cite{Chambolle04}.

Of course, in practice, Chambolle's (or any other iterative) algorithm
can only run for a finite number of iterations, thus the minimization
in line 11 of PIDAL-TV can only be solved approximately. However,
as stated in Theorem~\ref{th:Eckstein}, this will not compromise the
convergence of ADMM/PIDAL-TV, if the corresponding error sequence is
summable. To achieve this goal, we adopt a simple procedure in which
the internal variables of Chambolle's algorithm (the discrete gradient, see
\cite{Chambolle04}) are initialized, in each iteration of PIDAL-TV,
with those obtained in the previous iteration. We will now formalize
this idea and provide experimental evidence that this procedure does
produce a summable error sequence.

Let us define $\beta = \tau/\mu$ and let $({\bf s},{\bf q}) = \widetilde{\bPsi}^{(t)}_{\beta\, \mbox{\scriptsize TV}}({\bf r},{\bf p})$ be the result of running $t$
iterations of Chambolle's algorithm with its internal variables initialized
at ${\bf p}$, where ${\bf s}$ is the obtained (denoised) image
(which is approximately $\bPsi_{\beta\, \mbox{\scriptsize TV}}({\bf r})$)
and  ${\bf q}$ the final values of the internal variables.
Consider now two possible implementations of line 11 of PIDAL-TV:
\begin{align}
11 \mbox{(a)}. &\phantom{aaaaa} ({\bf u}_{k+1}^{(2,a)},{\bf p}_{k+1}) & = &\phantom{cc}
\widetilde{\bPsi}^{(t)}_{\beta\, \mbox{\scriptsize TV}}(\bnu^{(2)}_k,{\bf p}_k)\label{eq:11a}\\
11 \mbox{(b)}. &\phantom{aaaaa} ({\bf u}_{k+1}^{(2,b)}, \cdot) & =  & \phantom{cc}
\widetilde{\bPsi}^{(t)}_{\beta\, \mbox{\scriptsize TV}}(\bnu^{(2)}_k,{\bf p}).
\end{align}
Implementation $11$(a) uses the proposed internal variables initialization,
whereas in $11$(b) the internal variables are always initialized
to the same values (usually zeros).
Consider now the corresponding error sequences
\begin{align}
\rho^{(a)}_k & = \left\| {\bf u}_{k+1}^{(2,a)}
- \bPsi_{\beta\, \mbox{\scriptsize TV}} (\bnu_k^{(2)}) \right\| \\
\rho^{(b)}_k & = \left\| {\bf u}_{k+1}^{(2,b)}
- \bPsi_{\beta\, \mbox{\scriptsize TV}} (\bnu_k^{(2)}) \right\| .
\end{align}
Notice that since the two other minimizations (lines 9 and 13)
are solved exactly, the sequences $\rho^{(a)}_k$ and $\rho^{(b)}_k$
correspond to the sequence $\rho_k$ in Theorem~\ref{th:Eckstein}.

\begin{figure}
\centerline{\includegraphics[width=0.9\columnwidth]{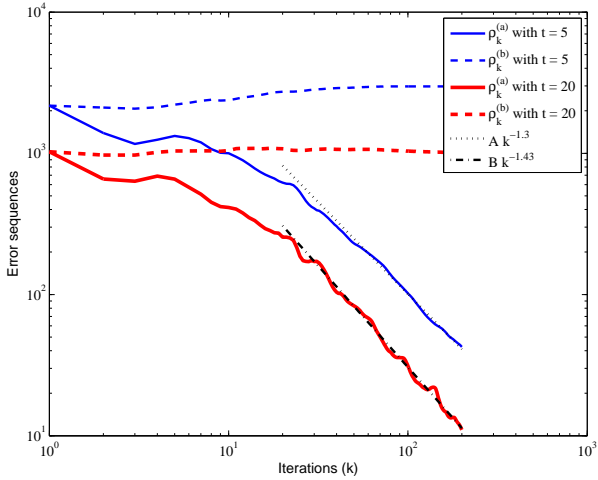}}
\caption{Error sequences $\rho_k^{(a)}$ and $\rho_k^{(b)}$ for $t=5$ and
$t=20$ (number of iterations of Chambolle's algorithm) and fitted functions
of the form $A\, t^\omega$ to the sequences $\{\rho_k^{(a)}, k=20,21,...200\}$.}
\label{fig:chambolle_init_errors}
\end{figure}

The following experiment provides evidence that $\rho^{(a)}_k$
is summable, but $\rho^{(b)}_k$ is not. Consider the same
setup as in the first experiment in~\cite{SetzerSteidlTeuber}:
the original image is a portion of the Cameraman image
scaled to a maximum value of 3000 and then blurred with
a Gaussian kernel of unit variance; the observed image is generated
according to (\ref{Poisson}). As in~\cite{SetzerSteidlTeuber}, we
set  $\tau = 0.008$ and  $\mu = \tau/50$. The number of iterations of Chambolle's
algorithm is set to 5 or 20. To compute
$\bPsi_{\beta\, \mbox{\scriptsize TV}} (\bnu_k^{(2)}) $ (almost)
exactly, we run 4000 iterations of Chambolle's algorithm.
In Figure~\ref{fig:chambolle_init_errors}, it is clear that the
$\rho^{(b)}_k$ sequences are not even decreasing, let alone
summable. In contrast, the sequences $\rho^{(a)}_k$ approach
zero, for both choices of $t$. Evidence for the summability
of the $\rho^{(a)}_k$ sequences is provided by the fact that
by fitting a function of the form $A\, (1/k)^\omega$ to the
tails of these sequences ({\it i.e.}, for $k=20,...,200$),
we obtain values of $\omega$ that are larger
than one ($\omega \simeq 1.3$, for $t=5$, and $\omega \simeq 1.43$,
for $t=20$).

In conclusion, the experiment reported in the previous paragraph,
though of course not a formal proof, strongly suggests that
by implementing line 11 of PIDAL-TV as in (\ref{eq:11a}), the
corresponding error sequence (with respect to the exact
minimizations) is summable, thus we can invoke Theorem~\ref{th:Eckstein}
to state that PIDAL-TV converges. Moreover, this experiment
shows that this is achieved with a quite small
number of iterations in each call of Chambolle's algorithm.
In all our experiments with PIDAL-TV, we thus use (\ref{eq:11a}),
with $t=5$.

\section{Poissonian Image Reconstruction with Frame-Based Regularization}
\label{sec:PIDAL-F}
We now consider the frame-based analysis criterion (\ref{eq:L_uncons_analysis}),
and the frame-based synthesis criterion (\ref{eq:L_uncons_synthesis}).

\subsection{Analysis Criterion}
In this case, the objective function is given by (\ref{eq:frame_analysis}),
which has the form (\ref{unconstrained_compound}) with  $J=3$,
\begin{equation}
g_1 \equiv {\cal L}, \;\;\; g_2 \equiv \tau  \|\cdot \|_1, \;\;\; g_3 \equiv \iota_{\mathbb{R}_+^n}
\end{equation}
and
\[
{\bf H}^{(1)} \equiv {\bf K}, \;\;\; {\bf H}^{(2)} \equiv {\bf P}, \;\;\; {\bf H}^{(3)} \equiv {\bf I}.
\]
The resulting instance of ADMM, which we call PIDAL-FA (where FA stands for ``frame analysis"),
is shown in Fig.~\ref{fig:PIDAL-FA}. The matrix being inverted in line 7 results from
assuming that  ${\bf P}$ is the analysis operator of a 1-tight (Parseval) frame, thus
 $({\bf H}^{(2)})^T{\bf H}^{(2)} = {\bf P}^T{\bf P} = {\bf I}$.
Notice that line 7 of PIDAL-FA corresponds to (\ref{eq:quadratic_problem})
for the particular form of matrix  ${\bf G}$ in this case,
${\bf G} = [{\bf K}^T \;\; {\bf P}^T\;\; {\bf I}\, ]^T$, which of course has
full column rank. As in PIDAL-TV, if $\bf K$ models a convolution,
the inverse $({\bf K}^T {\bf K} + 2\,{\bf I})^{-1}$ can be computed
with $O(n\log n)$ cost, using the FFT  (see (\ref{eq:Wiener2})).

\begin{figure}
\begin{center}
\fbox{\parbox{0.9\columnwidth}{\footnotesize
\begin{algorithm}{Poisson Image Deconvolution by AL (PIDAL-FA)}{}
Choose $\bo{u}_0^{(1)}$, $\bo{u}_0^{(2)}$, $\bo{u}_0^{(3)}$, $\bo{d}^{(1)}_{0}$, $\bo{d}^{(2)}_{0}$,
$\bo{d}^{(3)}_{0}$, $\mu$, and $\tau$. Set $k \leftarrow 0$.\\
\qrepeat \\
  $\bzeta_k^{(1)} \leftarrow {\bf u}_k^{(1)} + \bo{d}^{(1)}_{k}$\rule[-0.3cm]{0cm}{0.6cm}\\
  $\bzeta_k^{(2)} \leftarrow {\bf u}_k^{(2)} + \bo{d}^{(2)}_{k}$\rule[-0.3cm]{0cm}{0.6cm}\\
  $\bzeta_k^{(3)} \leftarrow {\bf u}_k^{(3)} + \bo{d}^{(3)}_{k}$\rule[-0.3cm]{0cm}{0.6cm}\\
  $\bgamma_k \leftarrow  {\bf K}^T \bzeta_k^{(1)} + {\bf P}^T\bzeta_k^{(2)} + \bzeta_k^{(3)} $\rule[-0.3cm]{0cm}{0.6cm}\\
  $\bo{z}_{k+1} \leftarrow \left( {\bf K}^T{\bf K} + 2\, {\bf I}\right)^{-1} \bgamma_k \rule[-0.3cm]{0cm}{0.6cm}$\\
        $\bnu_k^{(1)} \leftarrow {\bf K \, z}_{k+1} - \bo{d}^{(1)}_{k}$\rule[-0.3cm]{0cm}{0.6cm}\\
${\displaystyle \bo{u}^{(1)}_{k+1} \leftarrow
\frac{1}{2} \left( \nu_{i,k}^{(1)} - \frac{1}{\mu} +
\sqrt{\left( \nu_{i,k}^{(1)}-\frac{1}{\mu} \right)^2 + 4\, y_i}\; \right)  }$\rule[-0.3cm]{0cm}{0.6cm}\\
   $\bnu_k^{(2)} \leftarrow {\bf P} {\bf z}_{k+1}-\bo{d}^{(2)}_{k}$\rule[-0.3cm]{0cm}{0.6cm}\\
${\displaystyle\bo{u}_{k+1}^{(2)} \leftarrow \arg\min_{\bo{v}} \frac{\mu}{2}\|
\bo{v}-\bnu_k^{(2)} \|^2 + \tau \,\|{\bf v}\|_1}$.\rule[-0.3cm]{0cm}{0.6cm}\\
$\bnu_k^{(3)} \leftarrow {\bf z}_{k+1} - \bo{d}^{(3)}_{k}$\rule[-0.3cm]{0cm}{0.6cm}\\
${\displaystyle\bo{u}_{k+1}^{(3)} \leftarrow \max \left\{ \bnu_{k}^{(3)} , 0 \right\} }$.\rule[-0.3cm]{0cm}{0.6cm}\\
$\bo{d}^{(1)}_{k+1}  \leftarrow \bo{d}^{(1)}_{k} - \bo{K\, z}_{k+1} +\bo{u}_{k+1}^{(1)}$\rule[-0.3cm]{0cm}{0.6cm}\\
$\bo{d}^{(2)}_{k+1}   \leftarrow  \bo{d}^{(2)}_{k} - {\bf P} \bo{z}_{k+1} +\bo{u}_{k+1}^{(2)}$\rule[-0.3cm]{0cm}{0.6cm}\\
$\bo{d}^{(3)}_{k+1}   \leftarrow  \bo{d}^{(3)}_{k} -\bo{z}_{k+1} +\bo{u}_{k+1}^{(3)}$\rule[-0.3cm]{0cm}{0.6cm}\\
$k := k + 1$
\quntil some stopping criterion is satisfied.
 \end{algorithm} }}
 \caption{The PIDAL-FA algorithm.}\label{fig:PIDAL-FA}
\end{center}
\end{figure}

For most tight frames used in image processing,  products by ${\bf P}$
and ${\bf P}^T$ correspond to the inverse and direct transforms for
which fast algorithms exist. For example, in the case of
translation-invariant wavelet transforms, these products can
be computed using the undecimated wavelet transform with
$O(n\log n)$ cost \cite{lang}, \cite{Mallat}.
Curvelets also constitute a Parseval frame for which fast
$O(n\log n)$ implementations of the forward and inverse transform exist
\cite{CandesDemanetDonohoYing}. Yet another example of a redundant Parseval
frame is provided by complex wavelets, with the corresponding direct and inverse
transforms having $O(n)$ cost \cite{Kingsbury}, \cite{Selesnick}.
In conclusion, for a large class of choices of ${\bf P}$, the cost of lines
6, 10, and 15 of PIDAL-FA is $O(n\log n)$.

The expressions in lines 9 and 13 of PIDAL-FA are similar to those in lines 9 and 13
of PIDAL-TV, respectively; see also (\ref{eq:Poisson_Moreau}), (\ref{eq:Poisson_Moreau2}),
and (\ref{eq:orthant}).

The minimization in line 11 is, by definition, the Moreau proximity operator
of the $\ell_1$ norm \cite{CombettesSIAM}, which corresponds to a soft-threshold
(\ref{eq:soft}).

In summary, from the observations in the previous paragraphs, the computational costs of
the lines of PIDAL-FA are  the following. Lines 3, 4, 5, 9, 11, 12, 13, and 16
have $O(n)$ cost. Lines 6, 7, 8, 10, 14, and 15 have $O(n\log n)$ cost. Thus the
computational cost of PIDAL-FA scales as $O(n\log n)$.

Finally, convergence of PIDAL-FA is addressed by the following corollary  of Theorem \ref{th:Eckstein}.

\vspace{0.3cm}
\begin{corollary}\label{cor:PIDAL_FA}
{ \sl
The PIDAL-FA algorithm converges to a minimizer of  (\ref{eq:L_uncons_analysis}), provided one exists.}
\end{corollary}

\vspace{0.3cm}
\begin{proof}
The proof is similar to, but simpler than, that of Corollary~\ref{cor:PIDAL_TV},
since all the minimizations involved are solved exactly in closed form. Clearly,
matrix ${\bf G} = [{\bf K}^T \;\; {\bf P}^T\;\; {\bf I}\, ]^T$ has full column rank,
thus Theorem \ref{th:Eckstein} guarantees convergence to a minimizer of the objective
function.
\end{proof}
\vspace{0.3cm}

\subsection{Synthesis Criterion}
In the synthesis formulation,
the objective function is given by (\ref{eq:L_uncons_synthesis}),
which has the form (\ref{unconstrained_compound}) with  $J=3$,
\[
g_1 \equiv {\cal L}, \;\;\; g_2 \equiv \tau  \|\cdot \|_1, \;\;\; g_3 \equiv \iota_{\mathbb{R}_+^n}
\]
and
\[
{\bf H}^{(1)} \equiv {\bf K\, W}, \;\;\; {\bf H}^{(2)} \equiv {\bf I}, \;\;\; {\bf H}^{(3)} \equiv {\bf W}.
\]
The resulting ADMM algorithm, which we call PIDAL-FS (where FS stands for ``frame synthesis"),
is shown in Fig.~\ref{fig:PIDAL-FS}.

\begin{figure}
\begin{center}
\fbox{\parbox{0.9\columnwidth}{\footnotesize
\begin{algorithm}{Poisson Image Deconvolution by AL (PIDAL-FS)}{}
Choose $\bo{u}_0^{(1)}$, $\bo{u}_0^{(2)}$, $\bo{u}_0^{(3)}$, $\bo{d}^{(1)}_{0}$, $\bo{d}^{(2)}_{0}$,
$\bo{d}^{(3)}_{0}$, $\mu$, and $\tau$. Set $k \leftarrow 0$.\\
\qrepeat \\
  $\bzeta_k^{(1)} \leftarrow {\bf u}_k^{(1)} + \bo{d}^{(1)}_{k}$\rule[-0.3cm]{0cm}{0.6cm}\\
  $\bzeta_k^{(2)} \leftarrow {\bf u}_k^{(2)} + \bo{d}^{(2)}_{k}$\rule[-0.3cm]{0cm}{0.6cm}\\
  $\bzeta_k^{(3)} \leftarrow {\bf u}_k^{(3)} + \bo{d}^{(3)}_{k}$\rule[-0.3cm]{0cm}{0.6cm}\\
  $\bgamma_k \leftarrow  {\bf W}^T {\bf K}^T \bzeta_k^{(1)} + \bzeta_k^{(2)} + {\bf W}^T\bzeta_k^{(3)} $\rule[-0.3cm]{0cm}{0.6cm}\\
  $\bo{z}_{k+1} \leftarrow \left( {\bf W}^T {\bf K}^T{\bf K\, W} +  {\bf I} + {\bf W}^T{\bf W}\right)^{-1} \bgamma_k \rule[-0.3cm]{0cm}{0.6cm}$\\
        $\bnu_k^{(1)} \leftarrow {\bf K\, W\, z}_{k+1} - \bo{d}^{(1)}_{k}$\rule[-0.3cm]{0cm}{0.6cm}\\
${\displaystyle \bo{u}^{(1)}_{k+1} \leftarrow \frac{1}{2} \left( \nu_{i,k}^{(1)} - \frac{1}{\mu} +
\sqrt{\left( \nu_{i,k}^{(1)}-\frac{1}{\mu} \right)^2 + 4\, y_i}\; \right) }$\rule[-0.3cm]{0cm}{0.6cm}\\
   $\bnu_k^{(2)} \leftarrow  {\bf z}_{k+1}-\bo{d}^{(2)}_{k}$\rule[-0.3cm]{0cm}{0.6cm}\\
$\bo{u}_{k+1}^{(2)} \leftarrow \mbox{soft}(\bnu_k^{(2)} , \tau/\mu) $ \rule[-0.3cm]{0cm}{0.6cm}\\
$\bnu_k^{(3)} \leftarrow {\bf W\, z}_{k+1} - \bo{d}^{(3)}_{k}$\rule[-0.3cm]{0cm}{0.6cm}\\
${\displaystyle\bo{u}_{k+1}^{(3)} \leftarrow \max \left\{ \bnu_{k}^{(3)} , 0 \right\} }$.\rule[-0.3cm]{0cm}{0.6cm}\\
$\bo{d}^{(1)}_{k+1}  \leftarrow \bo{d}^{(1)}_{k} - \bo{K\, W\, z}_{k+1} + \bo{u}_{k+1}^{(1)}$\rule[-0.3cm]{0cm}{0.6cm}\\
$\bo{d}^{(2)}_{k+1}   \leftarrow  \bo{d}^{(2)}_{k} - \bo{z}_{k+1} + \bo{u}_{k+1}^{(2)}$\rule[-0.3cm]{0cm}{0.6cm}\\
$\bo{d}^{(3)}_{k+1}   \leftarrow  \bo{d}^{(3)}_{k} - \bo{W\, z}_{k+1} + \bo{u}_{k+1}^{(3)}$\rule[-0.3cm]{0cm}{0.6cm}\\
$k := k + 1$
\quntil some stopping criterion is satisfied.
 \end{algorithm} }}
 \caption{The PIDAL-FS algorithm.}\label{fig:PIDAL-FS}
\end{center}
\end{figure}

Notice that line 7 of PIDAL-FS corresponds to (\ref{eq:quadratic_problem})
for the particular form of matrix  ${\bf G}$ in this problem:
${\bf G} = [({\bf K \, W})^T \;\, {\bf I}\;\, {\bf W}^T\, ]^T$.
This matrix has of course full column rank. However, even if $\bf K$ models a
periodic convolution (thus is block circulant), the question remains of
how to efficiently compute the matrix inverse in line 7,
since ${\bf K\, W}$ is not block circulant. The next paragraph shows
how to sidestep this difficulty.

Consider that matrix ${\bf W}$ corresponds to a 1-tight (Parseval)
frame, {\it i.e.}, ${\bf W\,W}^H = {\bf I}$, and start by
noticing that $ {\bf W}^T {\bf K}^T{\bf K\, W} +  {\bf I} + {\bf W}^T{\bf W} =
{\bf W}^H({\bf K}^H{\bf K} + {\bf I}){\bf W} + {\bf I}$. Applying the
Sherman–-Morrison-–Woodbury (SMW) matrix inversion formula   yields
\begin{eqnarray}
\lefteqn{\left(  {\bf W}^H({\bf K}^H{\bf K} + {\bf I}){\bf W} + {\bf I} \right)^{-1}  =}  \nonumber \\
& = &  {\bf I} -   {\bf W}^H  \Bigl( {\bf W\,W}^H  + \bigl( {\bf K}^H  {\bf K}
+ {\bf I} \bigr)^{-1} \Bigr)^{-1}  {\bf W} \nonumber \\
& = &  {\bf I} -   {\bf W}^H  \Bigl( {\bf I}  + \bigl( {\bf K}^H  {\bf K}
+ {\bf I} \bigr)^{-1} \Bigr)^{-1}  {\bf W} .\label{eq:wiener_salsa_fs}
\end{eqnarray}
Using the factorization (\ref{eq:diagonalization}), we have
\begin{equation}
\Bigl( {\bf I}  + \bigl( {\bf K}^H  {\bf K}
+ {\bf I} \bigr)^{-1} \Bigr)^{-1} =  {\bf U}^H \bigl( {\bf I} + (|{\bf D}|^2 + {\bf I})^{-1} \bigr)^{-1} {\bf U},
 \label{eq:fft_filter}
\end{equation}
where both inversions have  $O(n)$ cost since $|{\bf D}|^2$ and
${\bf I}$ are diagonal, thus products by the matrix in
\eqref{eq:fft_filter} have the $O(n\log n)$ cost associated
to the FFT implementation of the products by ${\bf U}$ and ${\bf U}^H$.

The leading cost of line 7 of PIDAL-FS (given by \eqref{eq:wiener_salsa_fs}) will thus be either $O(n\log n)$
or the cost of the products by ${\bf W}^H$ and ${\bf W}$. As mentioned above,
for a large class of choices of frames, matrix-vector products by ${\bf W}$ and ${\bf W}^H$
have $O(n\log n)$ cost.

From the observations in the previous paragraphs, the computational costs of
the lines of PIDAL-FS are  the following. Lines 3, 4, 5, 9, 10, 11, 12, 13, 15, and 16
have $O(n)$ cost. Lines 6, 7, 8,  and 14 have $O(n\log n)$ cost. Thus the
computational cost of PIDAL-FS scales as $O(n\log n)$.

Finally, convergence of PIDAL-FS is
addressed by the following corollary  of Theorem \ref{th:Eckstein}.

\vspace{0.3cm}
\begin{corollary}
{ \sl
The PIDAL-FS algorithm converges to a minimizer of  (\ref{eq:L_uncons_synthesis}), provided one exists.}
\end{corollary}

\vspace{0.3cm}
\begin{proof}
The proof is similar to that of Corollary~\ref{cor:PIDAL_FA},
since all the minimizations involved are solved exactly in closed form. Clearly,
matrix ${\bf G} = [{\bf W}^T{\bf K}^T \;\; {\bf I}\;\; {\bf W}^T\, ]^T$ has full column rank,
thus Theorem \ref{th:Eckstein} guarantees convergence to a minimizer of the objective
function.
\end{proof}
\vspace{0.3cm}

\section{Experiments}
\label{sec:experiments}
We now report experiments where  PIDAL is compared with other
state-of-the-art methods, namely those proposed in \cite{DupeFadiliStarck},
\cite{Foi_ICIP2005}, \cite{SetzerSteidlTeuber}. All the algorithms
are implemented in MATLAB and the experiments are carried out
on a PC with a 3.0GHz Intel Core2Extreme CPU, with 4Gb of RAM,
under Microsoft Windows Vista. Unless otherwise indicated, we
adjust the regularization parameter $\tau$ to achieve the highest
{\it improvement in signal-noise-ratio} (ISNR $=10\log_{10}
\left(\|{\bf y}-{\bf x}\|_2^2/\|\widehat{\bf x}-{\bf x}\|_2^2\right)$).
The PIDAL algorithms are initialized as shown in Table~\ref{tab:init}.

\begin{table}
\caption{Initialization of the PIDAL algorithms.}\label{tab:init}
\begin{center}
\begin{tabular}{c|c|c|c|c|c|c}
 & $\bo{u}_0^{(1)}$ & $\bo{u}_0^{(2)}$ & $\bo{u}_0^{(3)}$ & $\bo{d}^{(1)}_{0}$ &
$\bo{d}^{(2)}_{0}$ & $\bo{d}^{(3)}_{0}$ \\ \hline
PIDAL-TV & ${\bf y}$ & ${\bf y}$ & ${\bf y}$ & ${\bf 0}$& ${\bf 0}$& ${\bf 0}$\\ \hline
PIDAL-FA & ${\bf y}$ & ${\bf P y}$ & ${\bf y}$ & ${\bf 0}$& ${\bf 0}$& ${\bf 0}$\\ \hline
PIDAL-FS & ${\bf y}$ & ${\bf W}^T{\bf K}^T{\bf y}$ & ${\bf K}^T{\bf y}$ & ${\bf 0}$& ${\bf 0}$& ${\bf 0}$\\ \hline
\end{tabular}
\end{center}
\end{table}

According to Theorem~\ref{th:Eckstein}, ADMM (thus PIDAL) converges
for any choice of $\mu > 0$. However, this parameter does influence
the speed of the algorithms. To our knowledge, there is work on
methods to choose this parameter for optimal speed; in our experiments,
we use the following rule of thumb, found to achieve satisfactory
results: $\mu = 60\, \tau/M$, where $M$ is the maximum intensity of
the original image. We have observed that the results do not change
significantly if this parameter is changed to one order of
magnitude below or above this choice.

\subsection{Comparison with \cite{SetzerSteidlTeuber}}
\label{sec:steidl_exp}
\begin{figure}
\centerline{\includegraphics[width=0.8\columnwidth]{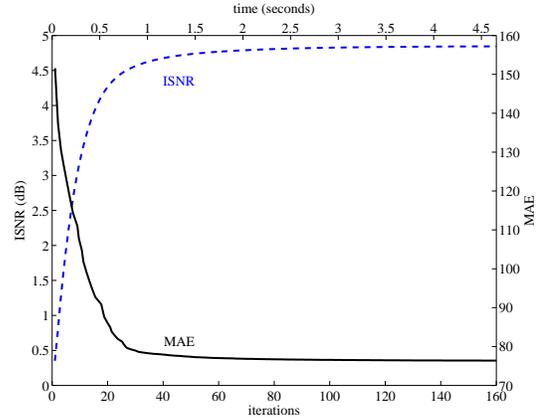}}
\caption{Evolution of the {\it mean absolute error} (MAE) and
{\it improvement in signal-noise-ratio} (ISNR) along the iterations
and elapsed time of PIDAL-TV, for the experiment of Section~\ref{sec:steidl_exp}. }
\label{fig:ISNR_MAE_Steidl}
\end{figure}

\begin{figure}
\centerline{\includegraphics[width=0.38\columnwidth]{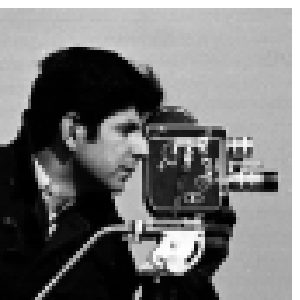}
\hspace{0.5mm}\includegraphics[width=0.38\columnwidth]{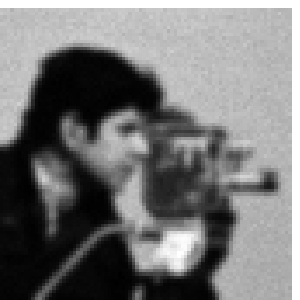}}

\vspace{1mm}
\centerline{\includegraphics[width=0.38\columnwidth]{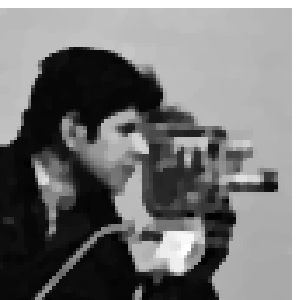}
\hspace{0.5mm}\includegraphics[width=0.38\columnwidth]{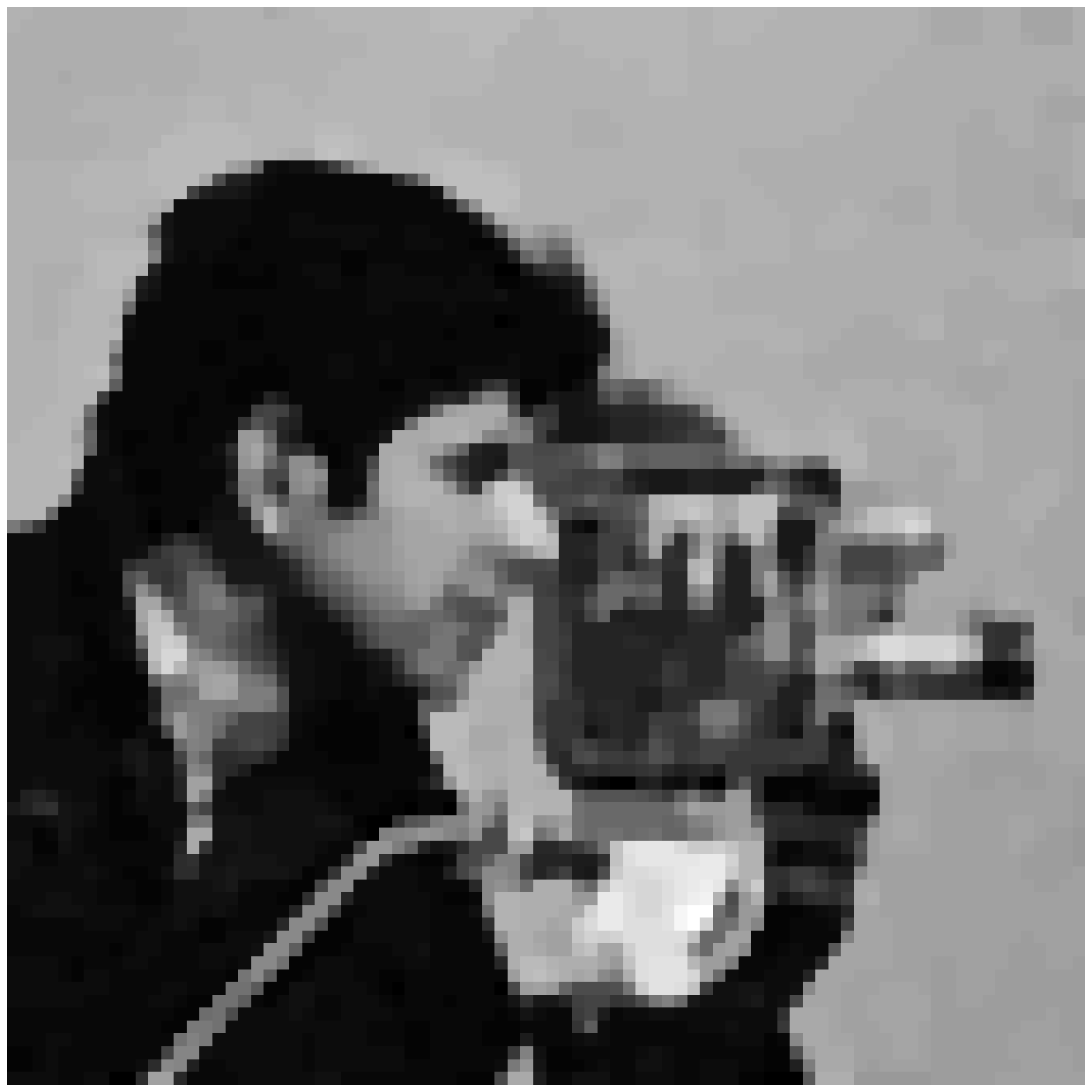}}
\caption{Experiment of Section~\ref{sec:steidl_exp}.
Top row: original (left); blurred and noisy image (right).
Bottom row: estimate from \cite{SetzerSteidlTeuber}; estimate by PIDAL-TV (ISNR=4.8dB).}
\label{fig:ISNR_MAE_Steidl_images_TV}
\end{figure}

We begin by comparing with the algorithms (PIDSplit and PIDSplit+)
proposed in \cite{SetzerSteidlTeuber}, which (as acknowledged by
the authors of \cite{SetzerSteidlTeuber}) is based on the earlier
version of PIDAL-TV \cite{FigueiredoDiasSSP09}.
The setup was already described in Section~\ref{sec:PIDAL-TV-conv-approx}:
the original image is a portion ($84\times 84$) of the Cameraman image, scaled to a
maximum value of 3000 and blurred with a Gaussian kernel of unit variance;
the observed image is generated according to (\ref{Poisson}); the
regularization parameter is set to $\tau = 0.008$.
In the experiments reported in
\cite{SetzerSteidlTeuber}, the TV denoising step of PIDAL-TV
is implemented by an inner iterative algorithm with a tight
stopping criterion based on the change between two consecutive
images. Our implementation of PIDAL-TV, as explained in
Section~\ref{sec:PIDAL-TV-conv-approx}, uses a small and
fixed number of iterations (just 5) of Chambolle's algorithm,
which is initialized as explained in that section. Because PIDSplit
and PIDSplit+ have no inner loop, each of its iterations
is roughly equivalent to one iteration of PIDAL-TV
with just one iteration of Chambolle's algorithm. In
\cite{SetzerSteidlTeuber}, PIDSplit and PIDSplit+ were run for
2150 iterations; we thus run PIDAL-TV for $2150/5 = 430$ iterations,
corresponding to roughly the same amount of computation.
Fig.~\ref{fig:ISNR_MAE_Steidl} shows the evolution of the
{\it mean absolute error} (MAE $= \|\widehat{\bf x}-{\bf x}\|_1/n$) and
ISNR along the first 160 iterations of PIDAL-TV (as well as elapsed time); it
is clear that convergence is achieved after less than 140
iterations (4.3 seconds, in our computer). This is dramatically
less than what is reported in \cite{SetzerSteidlTeuber} for
PIDAL-TV; in terms of iterations of PIDSplit and PIDSplit+,
this corresponds to approximately $150*5 = 750$ iterations,
thus also much less than the 2150 iterations (11 seconds)
reported in that work. Finally, Fig.~\ref{fig:ISNR_MAE_Steidl_images_TV}
shows the original, observed, and restored images; as expected,
the image estimates produced by  PIDSplit and PIDAL-TV are very similar.

\begin{figure}
\centerline{\includegraphics[width=0.8\columnwidth]{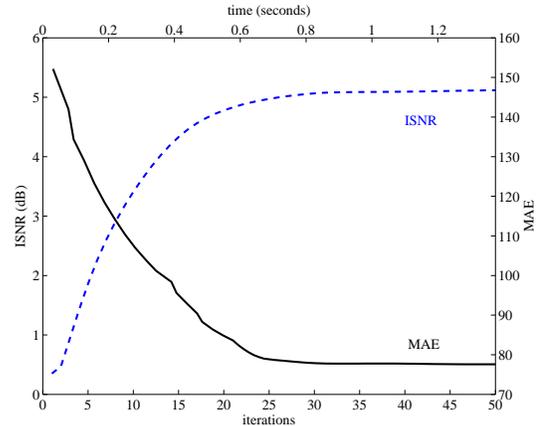}}
\caption{Experiment of Section~\ref{sec:steidl_exp}.
Evolution of the {\it mean absolute error} (MAE) and
{\it improvement in signal-noise-ratio} (ISNR) along the iterations
and elapsed time of PIDAL-FA.}
\label{fig:ISNR_MAE_Steidl_FA}
\end{figure}

\begin{figure}
\centerline{\includegraphics[width=0.8\columnwidth]{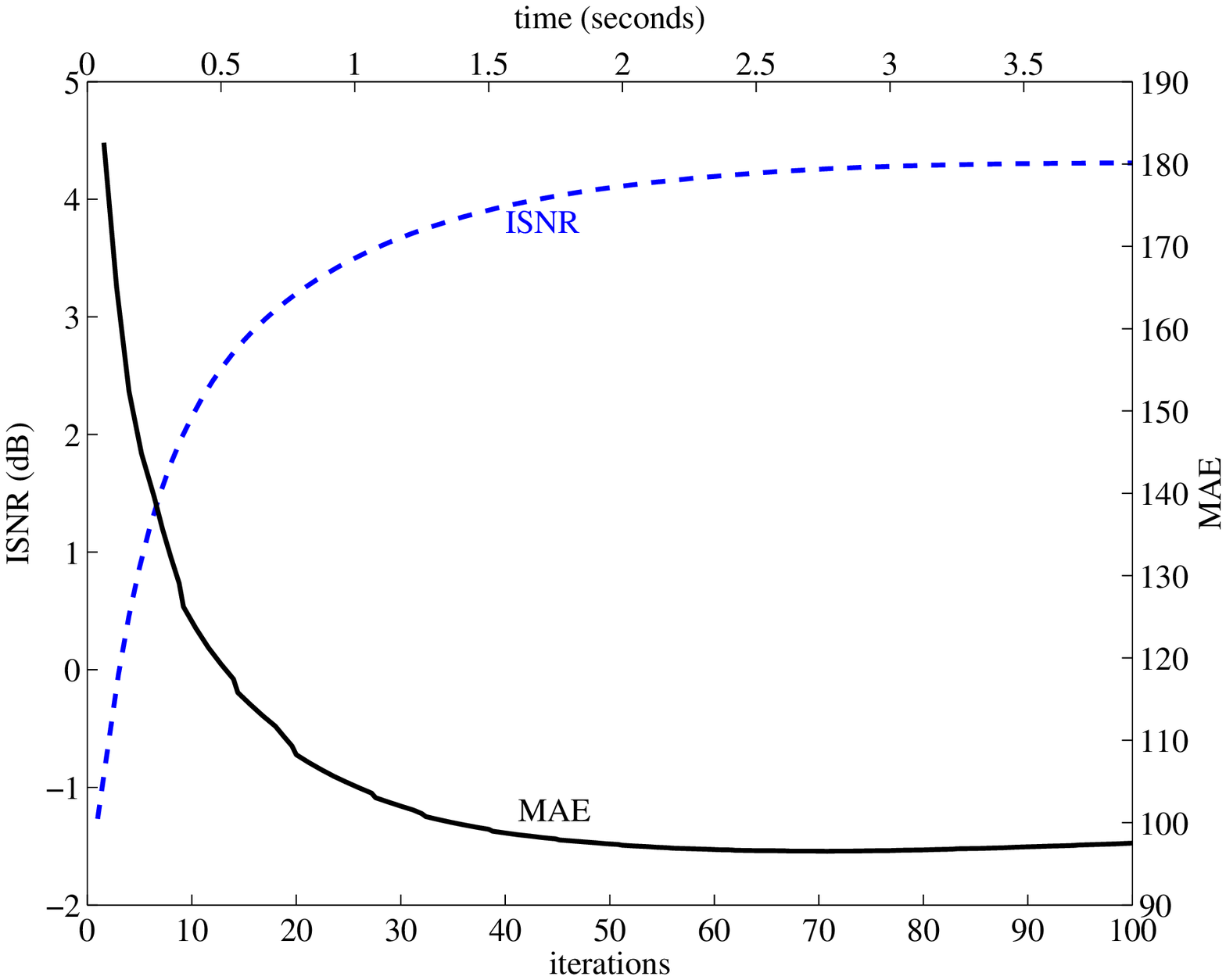}}
\caption{Experiment of Section~\ref{sec:steidl_exp}. Evolution of the {\it mean absolute error} (MAE) and
{\it improvement in signal-noise-ratio} (ISNR) along the iterations
and elapsed time of PIDAL-FS.}
\label{fig:ISNR_MAE_Steidl_FS}
\end{figure}

\begin{figure}
\centerline{\includegraphics[width=0.38\columnwidth]{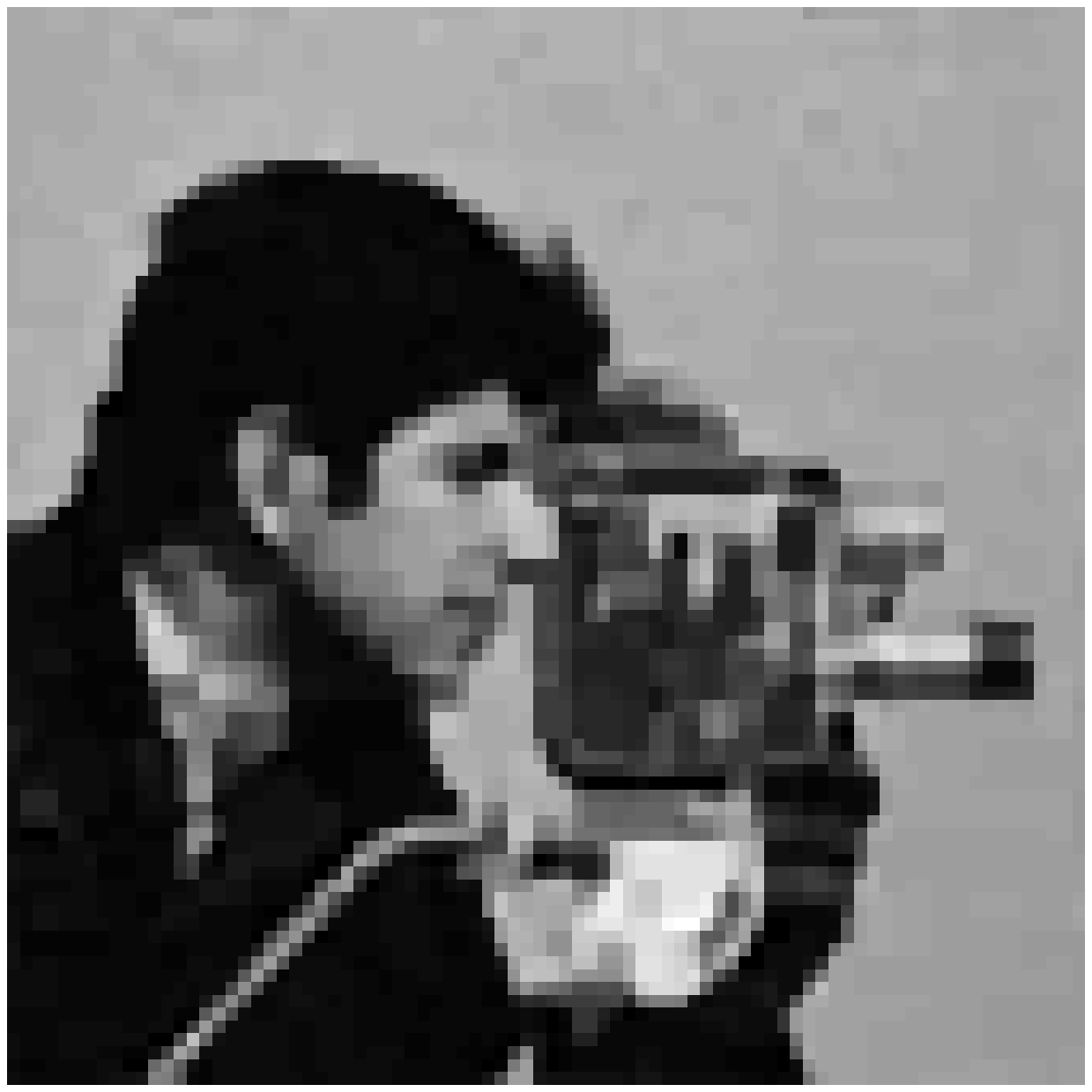}
\hspace{1mm}
\includegraphics[width=0.38\columnwidth]{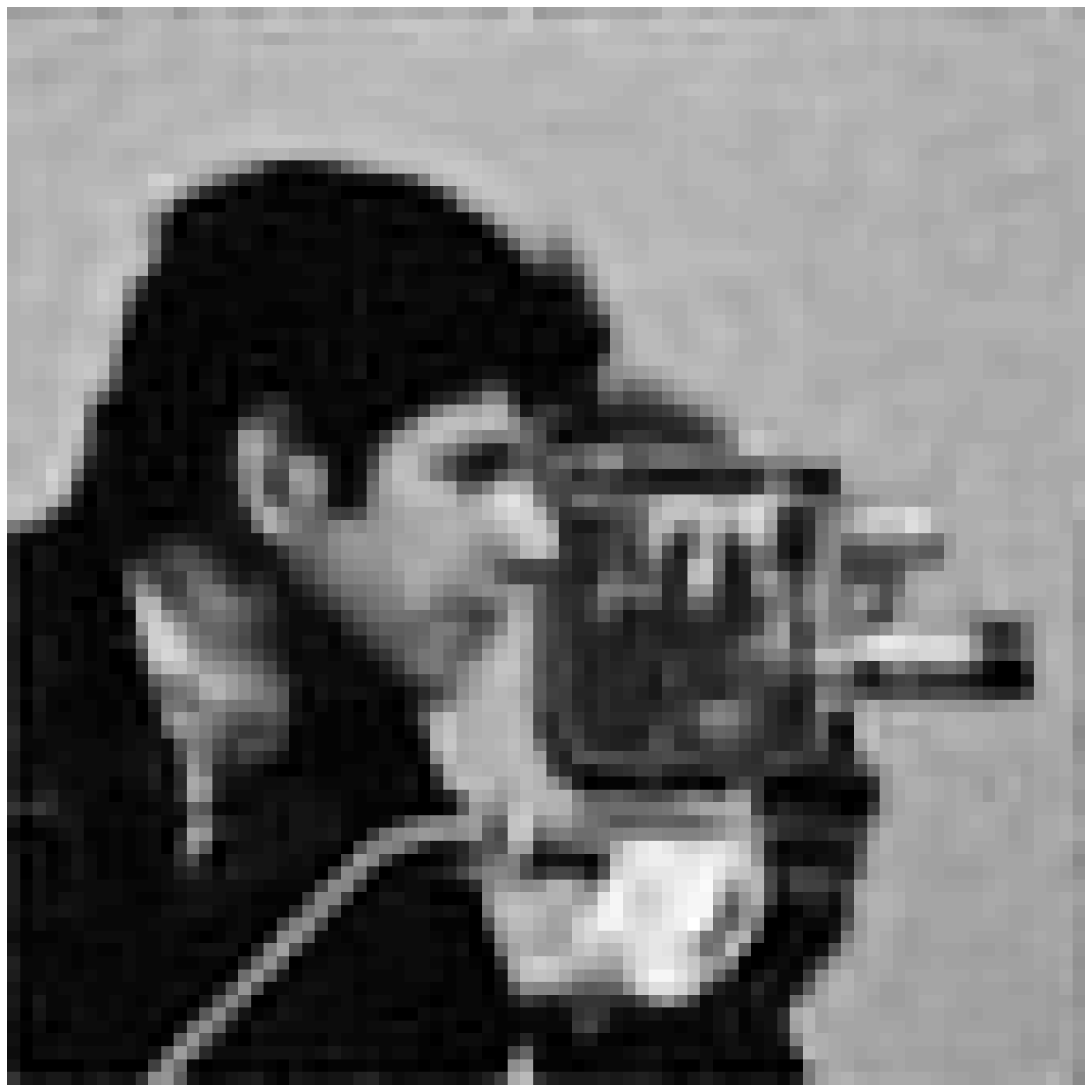}}
\caption{Experiment of Section~\ref{sec:steidl_exp}. Left: PIDAL-FA estimate (ISNR = 5.3dB). Right: PIDAL-FS estimate (ISNR = 4.3dB).}
\label{fig:ISNR_MAE_Steidl_images}
\end{figure}

Finally, we also tested PIDAL-FA and PIDAL-FS on the same example,
using a fully redundant Haar frame. The plots of ISNR and MAE are
presented in Figs.~\ref{fig:ISNR_MAE_Steidl_FA} and~\ref{fig:ISNR_MAE_Steidl_FS},
while the estimated images are shown in Fig.~\ref{fig:ISNR_MAE_Steidl_images}.

These results show that, in this example, PIDAL-FA performs slightly better than
PIDAL-TV in terms of ISNR and similarly in terms of MAE, with PIDAL-FA
achieving its best estimate faster than PIDAL-TV. The synthesis-based
criterion implemented by PIDAL-FS is a little worse in terms of both
ISNR and MAE, and PIDAL-FS also takes longer than PIDAL-FA
to achieve its best estimate. This poorer performance of the
synthesis formulation (in line with recent results in \cite{SelesnickFigueiredo})
was also found in all the experiments reported below,
so we will only present results from PIDAL-TV and PIDAL-FA.

\subsection{Comparison with \cite{Foi_ICIP2005}}
\label{sec:Foi_exp}
\begin{figure}
\centerline{\includegraphics[width=0.8\columnwidth]{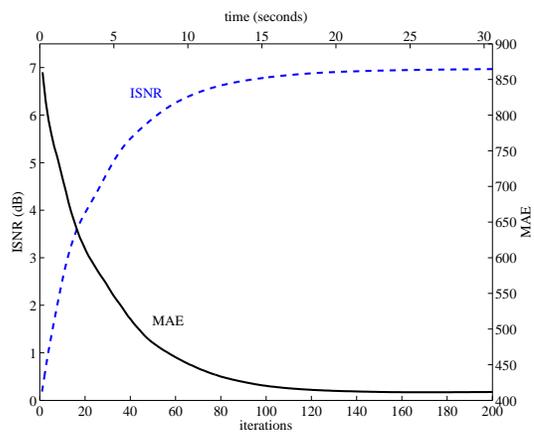}}
\caption{Experiment of Section~\ref{sec:Foi_exp}. Evolution of the {\it mean absolute error} (MAE) and
{\it improvement in signal-noise-ratio} (ISNR) along the iterations
and elapsed time of PIDAL-TV.}
\label{fig:ISNR_MAE_Foi_TV}
\end{figure}

\begin{figure}
\centerline{\includegraphics[width=0.8\columnwidth]{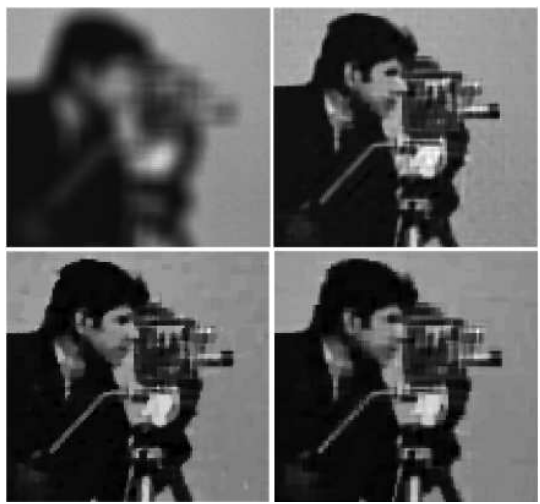}}
\caption{Experiment of Section~\ref{sec:Foi_exp}. Top row: blurred noisy image (left) and estimate from \cite{Foi_ICIP2005} (ISNR=6.61dB).
Bottom row: PIDAL-TV estimate (left, ISNR = 7.0dB); PIDAL-FA estimate (right, ISNR=6.95dB).}
\label{fig:Foi_images}
\end{figure}

The next experiment follows \cite{Foi_ICIP2005}: the original image
is the complete ($256\times 256$) Cameraman, scaled to a maximum value of 17600,
the blur is $9\times 9$ uniform. As in the experiment reported in the
previous subsection, this is a high SNR situation.
Fig.~\ref{fig:ISNR_MAE_Foi_TV} shows the evolution of the
MAE  and ISNR along the execution of PIDAL-TV; it
is clear that convergence is achieved after about 160
iterations (25 seconds, in our computer). A detail of the
blurred, and estimated images (from \cite{Foi_ICIP2005}
and using PIDAL-TV and PIDAL-FA) are shown in Fig.~\ref{fig:Foi_images}.
Although the TV and FA regularizers are considerably simpler than
the locally adaptive  approximation techniques used in \cite{Foi_ICIP2005},
both PIDAL-TV and PIDAL-FA achieve higher ISNR values (7.0dB and
6.95dB, respectively) than that reported in \cite{Foi_ICIP2005} (6.61dB).

\subsection{Comparison with \cite{DupeFadiliStarck}}
In the last set of experiments we compare our approach
with another recent state-of-the-art algorithm (herein
referred to as DFS), proposed in \cite{DupeFadiliStarck},
for which the MATLAB implementation is available at
{\tt\small www.greyc.ensicaen.fr/$\sim$fdupe}.
That work includes comparisons with other methods, namely:
Richardson-Lucy with multi-resolution support wavelet
regularization (RL-MRS) \cite{Starck1995};
fast translation invariant tree-pruning reconstruction
(FTITPR) \cite{Willett04}; Richardson-Lucy with
total variation regularization (RL-TV) \cite{Dey}.
The results in \cite{DupeFadiliStarck} show that the
algorithm therein proposed generally achieves better
performance ({\it i.e.}, lower MAE) than the others,
except for one of the images (a microscopy cell image)
where RL-MRS outperforms DFS. For this reason, we
will report results comparing PIDAL-TV and PIDAL-FA
versus DFS and RL-MRS. For PIDAL-FA, we use a
redundant Haar frame for the Cameraman image and
Daubechies-4 for the other images.
As in \cite{DupeFadiliStarck}, the original images
are scaled to a maximum value $M$,  belonging to
$\{5,30, 100, 255\}$, and then blurred by a $7\times 7$
uniform filter.

The DFS algorithm does not include a stopping criterion,
with the results reported in \cite{DupeFadiliStarck} having
been obtained by running a fixed number (200) of iterations.
In order to compare the running times of PIDAL-TV, PIDAL-FA,
and DFS, we run DFS until the MAE decreases less than 0.01\%
between two consecutive iterations.
Our algorithms are stopped when the following condition is met:
\[
\frac{\|{\bf z}_k - {\bf z}_{k-1}\|_2}{\|{\bf z}_{k-1}\|_2} \leq \delta,
\]
with $\delta = 0.005$ if $M=5$ and $\delta = 0.001$ in all the
other cases. Notice that this favors DFS, since a stopping criterion
based on MAE is not applicable in practice due to the absence of the
original image.

The results reported in Table~\ref{tab:results} show that, in 9
out of the 12 experiments, either PIDAL-TV or PIDAL-FA achieves
the lowest MAE. Notice however, that the main goal of this paper
was not to introduce a new restoration criterion aiming at obtaining
the lowest possible MAE (or any other performance measure), but rather
to introduce algorithms to solve the optimization problems
resulting from variational formulations of Poissonian image
restoration. In terms of computational efficiency, PIDAL-TV
and PIDAL-FA are clearly faster than the DRS algorithm, except
in the very low SNR situations ($M=5$) for two of the images (Cameraman
and Cell).

\begin{table*}
\centering
\caption{Results of the comparison of PIDAL-TV and PIDAL-FA with the algorithms proposed
in \cite{DupeFadiliStarck} and \cite{Starck1995} (average
over 10 runs); the reported times are in seconds.}\label{tab:results}
\vspace{0.3cm}
\begin{tabular}{ c | c || c | c | c | c | c | c | c | c | c | c }
\multicolumn{2}{}{} & \multicolumn{3}{|c|}{ PIDAL-TV} & \multicolumn{3}{|c|}{ PIDAL-FA} & \multicolumn{3}{|c|}{ DFS \cite{DupeFadiliStarck}} & RL-MRS \cite{Starck1995} \\  \hline
Image & $M$ & MAE & iterations & time & MAE & iterations & time & MAE & iterations & time & MAE \\  \hline \hline
Cameraman & 5 &  0.27 & 120 & 22 & {\bf 0.26} & 70 & 13 & 0.35 & 6 & 4.5 & 0.37 \\ \hline
Cameraman & 30 &  1.29 & 51 & 9.1 & {\bf 1.22} & 39 & 7.4 & 1.47 & 98 & 75 & 2.06 \\ \hline
Cameraman& 100 & 3.99 & 33  &  6.0 & {\bf 3.63} & 36 & 6.8 & 4.31 & 426 & 318 & 5.58 \\ \hline
 Cameraman & 255 & 8.99 & 32 & 5.8 & {\bf 8.45} & 37 & 7.0 & 10.26 & 480 & 358 & 12.3 \\ \hline
Neuron & 5 & {\bf 0.17} & 117 & 3.6  & 0.18 & 66 & 2.9 & 0.19 & 6 & 3.9 & 0.19 \\ \hline
Neuron & 30 & {\bf 0.68} & 54 & 1.8 & 0.77 & 44 & 2.0 & 0.82 & 161 & 77 &  0.95 \\ \hline
Neuron & 100 & {\bf 1.75} & 43 & 1.4 & 2.04 & 41 & 1.8 & 2.32 & 427 & 199 &  2.88 \\ \hline
Neuron & 255 & 3.52 & 43 & 1.4 & {\bf 3.47} & 42 & 1.9  & 5.25 & 202 & 97 &  6.31 \\ \hline
Cell & 5   & 0.12 & 56 & 10 & {\bf 0.11} & 36 & 7.6  & 0.12 & 6 & 4.5 &  0.12 \\ \hline
Cell & 30  & 0.57 & 31 & 6.5 & 0.54 & 39 & 8.2  & 0.56 & 85 & 64 &  {\bf 0.47} \\ \hline
Cell & 100 & 1.71 & 85 & 15 & 1.46 & 31 & 6.4  & 1.72 & 215 & 162 &  {\bf 1.37} \\ \hline
Cell & 255 & 3.77 & 89 & 17 & 3.33 & 34 & 7.0  & 5.45 & 410 & 308 &  {\bf 3.10} \\ \hline
\end{tabular}
\end{table*}

\section{Concluding Remarks}
We have propose new algorithms to handle the optimization problems
resulting from regularization approaches to the restoration of
Poissonian images. These optimization problems include several
difficulties: the Poisson log-likelihood is non-quadratic and
its gradient is not Lipschitz; the state-of-the-art regularizers
are non-smooth;
there is a non-negativity constraint. We have started by presenting
sufficient conditions for existence and uniqueness of solutions of
these optimization problems, for the following regularizers:
total-variation, frame-based analysis, and frame-based synthesis.
These problems were handled by adapting the alternating
direction method of multipliers (ADMM) to their particular forms.
 This adaptation is based on a new
way of using ADMM to deal with problems in which the objective
function is a linear combination of convex terms, which can be
used in many other problems. We gave sufficient conditions for
convergence and proved that these are met in the considered cases.
Finally, we have experimentally compared the proposed algorithms
against competing techniques, showing that they achieve state-of-the-art
performance both in terms of speed and restoration accuracy.

\section*{Appendix A: Convex Analysis}
We very briefly review some basic convex analysis results
used in this paper. For more details see \cite{CombettesSIAM,Zalinescu}.

Consider a function $f:{\cal X}\rightarrow \bar{\mathbb{R}} = \mathbb{R}\cup \{-\infty,+\infty\}$, where $\bar{\mathbb{R}}$ is called the extended real line, and ${\cal X}$ is a real Hilbert
space. The domain of function $f$ is the set $\mbox{dom}(f) = \{{\bf x}: f({\bf x})< + \infty\}$.

The function $f$ is {\it convex} if $f(\alpha{\bf u} + (1-\alpha){\bf v}) \leq
\alpha f({\bf u}) + (1-\alpha)f({\bf v})$, for any ${\bf u},{\bf v}\in
{\cal X}$ and any $\alpha\in\, [0,\, 1]$.  Convexity is said to be strict
if the inequality holds strictly ($<$) for any ${\bf u},{\bf v}\in
\mbox{dom}(f)$ and $\alpha\in\, ]0,\, 1[$.

The function is called {\it proper} if it is not equal to $+\infty$
everywhere and is never equal to  $-\infty$.

The function $f$ is {\it lower semi-continuous} (lsc) at ${\bf v}$ if
\[
\lim_{\delta\searrow 0}\; \inf_{{\bf x}\in B({\bf v},\delta)} f({\bf x}) \geq f({\bf v}),
\]
where $B({\bf v},\delta) = \{{\bf x}:\;\|{\bf x}-{\bf v}\| \leq \delta\}$
is the $\delta$-ball around ${\bf v}$, and $\|\cdot \|$ is
the norm in the Hilbert space ${\cal X}$. A function is called lsc if it is
lsc at every point of its domain.

A function $f$ is called {\it coercive} if it verifies
$\lim_{\|{\bf x}\| \rightarrow \infty} f({\bf x}) = +\infty$.
Proper, lsc, coercive functions play a key role in optimization via
the following theorem \cite{CombettesSIAM}:
\begin{theorem} If $f$ is a proper, lsc, coercive, convex function, then
 the set $\arg\min_{{\bf x}\in {\cal X}} f({\bf x})$ is nonempty.\label{th:propermin}
\end{theorem}

The next theorem concerns strictly convex functions.
\begin{theorem}
If $f$ is a strictly convex function, the set
$\arg\min_{{\bf x}\in {\cal X}} f({\bf x})$ possesses at most one element.\label{th:stritcconvex}
\end{theorem}

\end{document}